\documentclass[a4paper,11pt]{article}
\usepackage{amsmath}
\usepackage{amsfonts}
\usepackage{titlesec}
\usepackage{amssymb}
\usepackage{amsthm}
\usepackage{accents}

\usepackage{tikz,pgfplots}

\usepackage{geometry}
\usepackage{color}

\definecolor{pantone_312}{rgb}{0.0,0.6784,0.9778}
\definecolor{pantone_312_90}{rgb}{0.0,0.698,0.8745}
\definecolor{pantone_312_80}{rgb}{0.0,0.7216,0.8863}
\definecolor{pantone_312_70}{rgb}{0.0,0.749,0.8941}
\definecolor{pantone_312_60}{rgb}{0.2941,0.7804,0.9098}
\definecolor{pantone_312_50}{rgb}{0.4431,0.8078,0.9176}
\definecolor{pantone_312_40}{rgb}{0.5647,0.84310,0.9333}
\definecolor{pantone_312_30}{rgb}{0.6745,0.87840,0.9451}
\definecolor{pantone_312_20}{rgb}{0.7843,0.91760,0.9601}
\definecolor{pantone_312_10}{rgb}{0.8824,0.95290,0.9765}

\definecolor{pantone_315}{rgb}{0.0, 0.4314, 0.5373}

\definecolor{pantone_315_1}{rgb}{0.0, 0.2014, 0.3373}
\definecolor{pantone_315_2}{rgb}{0.0, 0.3314, 0.4373}
\definecolor{pantone_315_3}{rgb}{0.0, 0.3314, 0.5373}
\definecolor{pantone_315_4}{rgb}{0.0, 0.4314, 0.5373}
\definecolor{pantone_315_5}{rgb}{0.0, 0.5514, 0.6573}
\definecolor{pantone_315_6}{rgb}{0.0, 0.6514, 0.7573}
\definecolor{pantone_315_7}{rgb}{0.0, 0.6314, 0.8373}
\definecolor{pantone_315_8}{rgb}{0.0, 0.7314, 0.8373}
\definecolor{pantone_315_9}{rgb}{0.0, 0.8314, 0.9373}
\definecolor{pantone_315_10}{rgb}{0.0, 0.9314, 0.9373}

\definecolor{pantone_black_7}{rgb}{0.2431, 0.2431, 0.2314}

\definecolor{pantone_prozess_yellow}{rgb}{1.0, 0.9294, 0.0}
\definecolor{pantone_prozess_yellow_20}{rgb}{1.0, 0.9882, 0.8353}

\definecolor{pantone_369}{rgb}{0.4784, 0.7098, 0.0863}
\definecolor{pantone_369_50}{rgb}{0.75, 0.85, 0.55}
\definecolor{pantone_369_60}{rgb}{0.7, 0.8, 0.5}
\definecolor{pantone_369_10}{rgb}{0.9333, 0.9647, 0.9059}

\definecolor{light-gray}{gray}{0.69}
\definecolor{light-red}{rgb}{1.0,0.4,0.4}

\definecolor{light-blue}{rgb}{0.4,0.45,1}
\definecolor{light-green}{rgb}{0.5,0.8,0.0}
\definecolor{dark-green}{rgb}{0.0,0.4,0.0}
\definecolor{dark-red}{rgb}{1.0,0.3,0.3}

\newcommand{\DOmega}{U}

\newcommand*{\dt}[1]{\accentset{\mbox{\large\bfseries .}}{#1}}
\newcommand*{\ddt}[1]{\accentset{\mbox{\large\bfseries .\hspace{-0.25ex}.}}{#1}}

\usepackage[utf8]{inputenc}

\usepackage{mathtools}

\usepackage{array}
\makeatletter
\newcommand{\thickhline}{%
    \noalign {\ifnum 0=`}\fi \hrule height 1pt
    \futurelet \reserved@a \@xhline
}
\newcolumntype{"}{@{\hskip\tabcolsep\vrule width 1pt\hskip\tabcolsep}}
\makeatother

\newcommand{\ci}{\mathrm{i}} %% complex number "i"

\usepackage{framed}
\usepackage{graphicx, caption, subcaption}
\usepackage{hyperref}
\usepackage{dsfont}
\usepackage{float}

\usepackage{mhchem}

\newcommand{\quotes}[1]{``#1''}

\theoremstyle{definition}

\newtheorem{theorem}{Theorem}[section]

\newtheorem{definition}[theorem]{Definition}

\newcommand{\eps}{\varepsilon}

\newcommand{\rhophi}{\rho}%{\rho_{\phi}} %

\usepackage{diagbox}

\renewcommand{\thefootnote}{\fnsymbol{footnote}}
\renewcommand{\thefootnote}{\arabic{footnote}}

\begin{document}

\begin{center}
{\LARGE Shadow Lagrangian dynamics for superfluidity
\renewcommand{\thefootnote}{\fnsymbol{footnote}}\setcounter{footnote}{0}
 \hspace{-3pt}\footnote{P. Henning acknowledge the support by the Swedish Research Council (grant 2016-03339) and the G\"oran Gustafsson foundation and A.\ M.\ N. Niklasson is most grateful for the hospitality and pleasant environment at the Division 
of Scientific Computing at the department of Information Technology at Uppsala University, where he stayed during his 
participation of the development of this work. 
This work is supported by the U.S. Department of Energy, Office of Basic Energy Sciences (FWP LANLE8AN)
and by the U.S. Department of Energy through the Los Alamos National Laboratory.
Los Alamos National Laboratory is operated by Triad National Security, LLC, for the National Nuclear Security
Administration of the U.S. Department of Energy Contract No. 892333218NCA000001. The manuscript has the LA-UR number 19-32663.}}\\[2em]
\end{center}

\begin{center}
{\large Patrick Henning\footnote[1]{Faculty of Mathematics, Ruhr-University Bochum, DE-44801 Bochum, Germany and Department of Mathematics, KTH Royal Institute of Technology, SE-100 44 Stockholm, Sweden.} and 
Anders M. N. Niklasson\footnote[2]{Theoretical Division, Los Alamos National Laboratory, Los Alamos, New Mexico 87545, USA 
and Division of Scientific Computing, Department of Information Technology, Uppsala University Box 337, SE-751 05 Uppsala, Sweden
}}\\[2em]
\end{center}

%35Q55, 65N12, 65N25, 65N30, 81Q05

\begin{center}
{\large{\today}}
\end{center}

\begin{center}
\end{center}

\begin{abstract}
Motivated by a similar approach for Born-Oppenheimer molecular dynamics, this paper proposes an extended ``shadow'' Lagrangian density for quantum states of superfluids. The extended Lagrangian contains an additional field variable that is forced to follow the wave function of the quantum state through a rapidly oscillating extended harmonic oscillator. By considering the adiabatic limit for large frequencies of the harmonic oscillator, we can derive the two equations of motions, a Schr\"odinger-type equation for the quantum state and a wave equation for the extended field variable. The equations are coupled in a nonlinear way, but each equation individually is linear with respect to the variable that it defines. The computational advantage of this new system is that it can be easily discretized using linear time stepping methods, where we propose to use a Crank-Nicolson-type approach for the Schr\"odinger equation and an extended leapfrog scheme for the wave equation. Furthermore, the difference between the quantum state and the extended field variable defines a consistency error that should go to zero if the frequency tends to infinity. By coupling the time-step size in our discretization to the frequency of the harmonic oscillator we can extract an easily computable consistency error indicator that can be used to estimate the numerical error without any additional costs. The findings are illustrated in numerical experiments.
\end{abstract}

\section{Introduction}
The phenomenon of superfluidity is often, though not exclusively, studied and described through Bose-Einstein condensation \cite{ARV01,Fet09,MAH99}. A Bose-Einstein condensate (BEC) is a remarkable state of matter that is formed when a dilute gas of bosonic particles is cooled down to extremely low temperatures close to the absolute zero \cite{Bos24,DGP99,Ein24,PiS03}. In this case, most of the bosons occupy the same quantum state (the ground state) and loose their identity. In fact, the individual packages interfere and behave like a single macroscopic \quotes{super-particle}. A typical example for a superfluid Bose-Einstein condensate is Helium-4 (\ce{{}^4He}), which allows to study the arising quantum phenomena (such as vortices with a quantized circulation) on an observable scale. Mathematically, the behavior of Bose-Einstein condensates can be described by nonlinear Schr\"odinger equations \cite{Aft06,LSY01,PiS03}, were the Gross-Piatevskii equation (GPE) is by far the most common model. In non-dimensional form, the basic GPE reads 
\begin{equation}\label{NSE}
\ci \dt{\psi} = \tfrac{1}{2} \nabla^2 \psi + V\psi + \kappa |\psi|^2\psi,
\end{equation}
where $\psi(\mathbf{x},t)$ describes the time-evolution of the quantum state, $V(\mathbf{x})$ is a magnetic trapping potential and $\kappa$ is a parameter that characterizes particle interactions. Nonlinear Schr\"odinger equations such as the GPE \eqref{NSE} come with important invariants, such as the energy or the number of particles (cf.~\cite{BaC13b}). This needs to be considered in numerical computations by selecting a suitable time-discretization that avoids significant energy dissipation or energy blow-up over time. Common choices for time integrators that allow only for small (or no) oscillations of the energy are either symplectic methods (cf. \cite{HeM17,KaM99,San88} and the references therein) or schemes that conserve the initial discrete energy up to machine precision (cf. \cite{ABB13d,ADK91,BaC13,Bes04,HeP17,San84}). Such time integrators are typically implicit or semi-implicit which involves (except for the Besse relaxation scheme \cite{Bes04}) the solving of a non-linear problem in each time step, since the time-dependent variable ${\psi }$ appears in the last potential term as $|\psi|^2$. Using more simple time integrators that would make the computations per time step cheaper typically results in instabilities due the aforementioned energy blow-up (cf. the experiments in \cite{HeW19}) or it leads to numerical overdamping expressed by a significant loss of energy that makes the numerical approximations useless (cf. the numerical experiments in \cite{HeP17}).

The problem is similar to a common problem in first-principles Born--Oppenheimer molecular
dynamics, where the forces acting on the atomic degrees of freedom are given only at
the fully relaxed electronic ground state, which is determined by a time-independent
non-linear Schr\"{o}dinger-like equation. 
If the equilibrated ground state solution for the electron density is not sufficiently accurate the calculated forces are not conservative and instabilities and energy drift may occur, cf. \cite{DRemler90,PPulay04,ANiklasson07}.
A recent solution to this problem 
is based on an extended Lagrangian approach to Born--Oppenheimer molecular dynamics (XL-BOMD),
where additional auxiliary electronic dynamical field
variables are included in an extended Lagrangian formalism \cite{ANiklasson08,ANiklasson09,GZheng11,ANiklasson14,ANiklasson17}.
The auxiliary electronic variables are constrained to follow closely the exact electronic ground state.
A linearization of the underlying energy functional around the auxiliary eletronic degrees of freedom
can then be used to avoid the non-linear optimization without any significant loss of accuracy.
The basic idea is a backward error analysis, where instead of calculating
approximate forces through an iterative non-linear optimization procedure for an underlying exact
potential energy surface, we can calculate exact forces in a fast and simple way, but for an underlying approximate
``shadow'' potential energy surface. In a classical adiabatic limit the dynamics of XL-BOMD
is described by two coupled equations of motion that each are linear with respect to their time-dependent variables.
In this way the non-linear ground state optimization is avoided and the forces remain conservative for
an underlying approximate ``shadow'' potential energy surface, which is a close approximation to the exact Born--Oppenheimer
potential energy surface.  Our idea here is to use the same basic approach applied to the time-dependent
non-linear Schr\"{o}dinger equation, equation \eqref{NSE}.  This can be achieved by reformulating the Lagrangian,
\begin{equation}\label{LagDens}
{\cal L }_{\mbox{\tiny std}}(\psi,\dt{\psi})  = \int_{\DOmega} \frac{\ci}{2} \left(\psi^*\dt{\psi} - \dt{\psi}^*\psi\right) - \frac{1}{2} \nabla \psi^* \cdot \nabla \psi
- V\psi^* \psi - \frac{1}{2} \kappa (\psi^* \psi)^2 \hspace{2pt}dx,
\end{equation}
corresponding to equation \eqref{NSE}, where $U\subset \mathbb{R}^d$ denotes the physical domain (for $d\in \{1,2,3\}$).
The Lagrangian can be modified to include extended auxiliary electronic degrees of freedom, represented by $\phi$.
The extended field variable, $\phi$, can then be forced to follow $\psi$ through an extended harmonic oscillator.
Afterwards, we can linearize the potential energy term, $(\psi^* \psi)^2$, in equation\ \eqref{LagDens} around the approximate electronic degrees of freedom, $\phi$, that evolves through the harmonic oscillator
centered around $\psi$. This can be accomplished through our definition of an extended ``shadow'' Lagrangian,
\begin{eqnarray*}
\lefteqn{ {\cal L}(\psi,\dt{\psi},\phi,\dt{\phi}) = \int_{\DOmega} \frac{\ci}{2} \left(\psi^*\dt{\psi} - \dt{\psi}^*\psi\right)
- \frac{1}{2} \nabla \psi^* \cdot \nabla \psi - V \psi^* \psi \hspace{2pt} dx} \\
&\enspace& \quad + \int_{\DOmega} - \frac{1}{2} \kappa \left[(2\psi^* -\phi^*)\phi \phi^*(2\psi - \phi)\right] + \frac{\mu}{2}\dt{\phi}^*\dt{\phi} - \frac{\mu \omega^2}{2} \left[(\psi^*-\phi^*)(\psi - \phi)\right] \hspace{2pt} dx.
\end{eqnarray*}
Here $\mu$ is a fictitious electron mass parameter and $\omega$ is the frequency of the extended harmonic oscillator
that is given by the last two terms of the Lagrangian.
In analogy with XL-BOMD \cite{ANiklasson14,ANiklasson17} we can now derive the Euler-Lagrange equations of motion
in an adiabatic limit, where we assume that $\mu \rightarrow 0$ and $\omega \rightarrow \infty$ such that
$\mu \omega = {\rm constant}$. In analogy with XL-BOMD we further assert that 
$|\psi-\phi| \sim {\cal O}(\omega^{-2})$, which can be validated a posteriori \cite{ANiklasson14}. 
In this limit we find that the equations of motion are given by
\begin{equation}\label{XLNSE}
\begin{array}{l}
\ci\dt{\psi} = \frac{1}{2} \nabla^2 \psi + V\psi + \kappa |\phi|^2(2\psi-\phi),\\
~~\\
\ddt{\phi} = \omega^2(\psi-\phi).
\end{array}
\end{equation}
The first equation governs the evolution of $\psi(t)$ and is linear in $\psi$, but non-linear in $\phi$, whereas
the second equation that describes the time evolution of $\phi(t)$ is linear both in $\psi$ and $\phi$.
As long as $\phi(t)$ stays close to $\psi(t)$, the quasi-linearized set of equations
given by equation \eqref{XLNSE} will closely approximate the time evolution of $\psi(t)$ 
in the original non-linear Schr\"{o}dinger equation, equation \eqref{NSE}. The computational advantage of \eqref{XLNSE} is that it can be integrated implicitly in time in such a way that each time step only requires the solving of one {\it linear} system of equations. To achieve this, we discretize the first equation of \eqref{NSE} with an implicit method of choice, which leads to a linear problem for $\psi^{n+1}$, where $\psi^n$ and $\phi^n$ are given from the previous time step and the implicit contribution $\phi^{n+1}$ is unknown. However, this missing contribution can be inferred from the linear second equation, i.e., $\ddt{\phi} = \omega^2(\psi-\phi)$, which can be solved cheaply using a suitable explicit time step. We no longer require an iterative non-linear optimization procedure and the properties of the arising numerical scheme for \eqref{XLNSE} are essentially dominated by the choice of the numerical method for the first equation (\quotes{the implicit part}). In this way it is possible to avoid numerical instabilities,
a broken time-reversal symmetry and energy drift. At the same time the computational cost is significantly reduced. Furthermore, the difference between $\phi$ and $\psi$ introduces a consistency error that should decay proportional to $\omega^{-2}$. This observation can be directly exploited when discretizing the equations \eqref{XLNSE} in time, since it yields a computable error indicator \quotes{for free}. This is a significant advantage compared to alternative time-discretizations of nonlinear Schr\"odinger equations, where the computation of error indicators typically requires considerable additional costs or a repeated solving for different step sizes until a converged state is observed.

The extended Lagrangian approach to first-principles electronic structure theory was pioneered by Car and Parrinello
in 1985 \cite{RCar85}. Their extended Lagrangian approach is different from XL-BOMD and our shadow Lagrangian density
formulation for the non-linear Schr\"{o}dinger equation introduced here, but there are several interesting
similarities \cite{ANiklasson17}.

$\\$
{\bf Outline:} The structure of this paper is as follows. In Section \ref{section-extended-lagrangian} we formulate the extended shadow Lagrangian in a more general way. Furthermore, we identify the conserved energy associated with the shadow Lagrangian, we derive the arising extended equations of motion, and we introduce a corresponding adiabatic limit approximation to which we refer as the shadow Lagrangian equations for superfluidity. In Section \ref{section:discretization} we propose and discuss a semi-explicit time discretization of the shadow Lagrangian equations, which is numerically validated in Section \ref{section:numerical-experiments} using an additional finite element discretization in space.

\section{Extended Lagrangian}
\label{section-extended-lagrangian}
In this section we introduce an extended Lagrangian formulation for superfluidity. The setting will be more general than what we sketched in the introduction. This is to emphasize that the approach also can be applied to more complex models. After introducing the extended Lagrangian, we derive the governing partial differential equations and pass to the adiabatic limit in order to obtain the final shadow Lagrangian equations that we propose as a starting point for new discretization schemes.

\subsection{Notation}
For a formal description of the Lagrangian we will use standard notation for Lebesgue and Sobolev spaces. In the following $\DOmega \subset \mathbb{R}^d$ denotes a connected and bounded subset of $\mathbb{R}^d$, our computational domain. On $\DOmega$ we let $L^2(\DOmega)$ denote the space of complex-valued and square-integrable functions. The corresponding $L^2$-norm of a function $v\in L^2(\DOmega)$ is then given by $\| v \|_{L^2(\DOmega)}= \left(\int_{\DOmega} |v(x)|^2 \hspace{2pt} dx \right)^{1/2}$. The space $H^1_0(\DOmega)$ denotes the Sobolev space of (complex-valued) weakly-differentiable functions, where all partial weak derivatives are in $L^2(\DOmega)$. On $H^1_0(\DOmega)$, the standard $H^1$-norm is defined as
$$
\| v \|_{H^1(\DOmega)} = \left(\int_{\DOmega} |v(x)|^2 + |\nabla v(x)|^2 \hspace{2pt} dx \right)^{1/2}.
$$
Finally, we denote the duality pairing of an element of a Hilbert space and an element from its dual space by $\langle \cdot , \cdot \rangle$.

\subsection{Extended governing equations}

The standard Lagrangian $\mathcal{L}_{\mbox{\tiny std}}$ used to describe (superfluid) Bose-Einstein condensates at ultra-low temperatures which occupy a domain $\DOmega\subset \mathbb{R}^d$ takes the form
\begin{eqnarray*}
\mathcal{L}_{\mbox{\tiny std}}(\psi, \dt \psi ) = \int_{\DOmega} \frac{\ci}{2} (\psi^{\ast} \dt{\psi} -  \dt{\psi}^{\ast} \psi ) - \frac{1}{2} \nabla \psi^{\ast} \cdot \nabla \psi - \frac{1}{2} V[\rho_{\psi}] \hspace{2pt}\rho_{\psi}  \hspace{2pt} dx,
\end{eqnarray*}
where $\rho_{\psi}=|\psi|^2$ denotes the density of the quantum state $\psi$. The function $V[\rho_{\psi}]$ describes a general dependency on this density, where we indirectly assume that $V$ is such that $\int_{\DOmega} V[\rho_{\psi}] \hspace{2pt}\rho_{\psi} \hspace{2pt} dx$ is finite for all possible quantum states $\psi \in H^1_0(\DOmega)$ and that $V[\rho_{\psi}]$ is real-Fr\'echet differentiable with respect to $\rho_{\psi}$. The most common choice for $V[\rho_{\psi}]$ in the context of superfluids is the Gross-Pitaevskii model \cite{PiS03}, where we have $V[\rho_{\psi}] = 2 V_0(x) + \kappa \hspace{2pt} \rho_{\psi}$ for an external trapping potential $V_0(x)$ and $\kappa$ being a parameter that characterizes interactions between particles. In particular, $\kappa$ is proportional to the scattering length of the bosons that form the superfluid (Bose-Einstein condensate) and $\kappa$ also depends on the mass and the number of particles (cf. \cite{BaC13b} for details on the scaling of the equation in non-dimensional form).

As described in the introduction for the Gross-Pitaevskii equation, we shall extend the standard Lagrangian density by an artificial harmonic oscillator. We will also linearize the potential term $ V[\rho_{\psi}] \hspace{2pt}\rho_{\psi}$ (with respect to its dependency on the quantum state $\psi$ in the arising governing equations, cf. equation \eqref{eqn-for-psi-full} below). For that, we introduce an extended field variable $\phi$ that is forced to follow the quantum state $\psi$ through the extended harmonic oscillator that is oscillating around $\psi$. This can be modeled by an extended (shadow) Lagrangian 
$$
\mathcal{L} : [H^1_0(\DOmega)]^4 \rightarrow \mathbb{R}
$$
that is given by
\begin{eqnarray*}
\lefteqn{ \mathcal{L}(\psi, \phi , \dt\psi , \dt\phi ) = \int_{\DOmega} \frac{\ci}{2} (\psi^{\ast} \dt{\psi} -  \dt{\psi}^{\ast} \psi ) - \frac{1}{2} \nabla \psi^{\ast} \cdot \nabla  \psi - \frac{1}{2} (V[\rhophi] - \rhophi V^{\prime}[\rhophi]) \hspace{2pt}\psi^{\ast} \psi  \hspace{2pt} dx}\\
&\enspace& \quad + \int_{\DOmega} - \frac{1}{2} \rhophi V^{\prime}[\rhophi] \hspace{2pt}(2 \psi^{\ast} - \phi^{\ast}) (2 \psi - \phi) + \frac{\mu}{2} \dt\phi^{\ast} \dt{\phi} - \frac{\mu \omega^2}{ 2} [(\phi^{\ast} - \psi^{\ast})(\psi - \phi)]  \hspace{2pt} dx,
\end{eqnarray*}
where $\rhophi:=|\phi|^2$ denotes the density of $\phi$ and the constant $\mu$ is a fictitious mass parameter. The first four terms in the extended Lagrangian density $\mathcal{L}$ model the dynamics of the superfluid, where we observe that for $\psi=\phi$, these four terms collapse to the standard Lagrangian $\mathcal{L}_{\mbox{\tiny std}}$. 
Hence, as long as $\phi$ stays close to $\psi$, the difference in the potential energy will therefore remain small. This is achieved by propagating $\phi$ though a harmonic oscillator that is centered around $\psi$. The last two terms in $\mathcal{L}$ model precisely this harmonic oscillator whose frequency is denoted by $\omega$. For large frequencies, this makes sure that $\psi$ stays close to $\phi$ (with $|\psi-\phi|\sim \omega^{-2}$) and hence that the potential energy surface induced by $\phi$ is close to the exact energy surface.

In order to derive the governing equations that arise from this extended Lagrangian, we can compute the partial (real-)Fr\'echet derivatives of the operator $\mathcal{L}$ and exploit the Euler-Lagrange formalism
 $\frac{\mbox{d}}{\mbox{d}t} \partial_{\dt{\psi}} \mathcal{L}(\psi,\phi,\dt{\psi} ,\dt{\phi}) - \partial_{\psi} \mathcal{L}(\psi,\phi,\dt{\psi} ,\dt{\phi} ) = 0$. With this we obtain the equation of motion for the quantum state $\psi$ (in the variational sense) as
\begin{eqnarray}
\label{eqn-for-psi-full}
\ci \dt{\psi} = - \frac{1}{2}\nabla^2 \psi + \frac{1}{2} (V[\rhophi] - \rhophi V^{\prime}[\rhophi]) \psi  + \rhophi V^{\prime}[\rhophi] (2 \psi - \phi) + \frac{\mu \omega^2}{2} ( \psi - \phi).
\end{eqnarray}
Analogously, we conclude from $\frac{\mbox{d}}{\mbox{d}t} \partial_{\dt{\phi}} \mathcal{L}(\psi,\phi,\dt{\psi} ,\dt{\phi}) - \partial_{\phi} \mathcal{L}(\psi,\phi,\dt{\psi} ,\dt{\phi} ) = 0$ the Euler-Lagrange equation for the extended field variable $\phi$ with
\begin{eqnarray}
\label{eqn-for-phi-full}
 \mu \ddt{\phi} =  - (V^{\prime}[\rhophi] + \rhophi V^{\prime\prime}[\rhophi]) |2 \psi - \phi|^2  \phi + \rhophi V^{\prime}[\rhophi] (2 \psi - \phi  ) + \rhophi  V^{\prime\prime}[\rhophi] \rho_{\psi} \phi + \mu \omega^2 (\psi - \phi).
\end{eqnarray}
Since $\psi$ and $\phi$ should, in the adiabatic limit, converge to the same quantum state $u$, it is reasonable to select the initial values as 
$$
\psi(0)=\phi(0)=u_0,
$$ 
where $u_0$ denotes some known initial state. 

For the system \eqref{eqn-for-psi-full}-\eqref{eqn-for-phi-full} to be well-posed, we still require a second initial value, $\dt{\phi} (0)$, for the extended field variable whose dynamics are described by \eqref{eqn-for-phi-full}. This initial value is typically not explicitly given and its choice is part of the extended model. In order to derive a reasonable condition for  $\dt{\phi} (0)$, we shall first identify the conserved quantity of the extended Lagrangian, which is done in the next subsection.

\subsection{Conserved energy}
In the next step, we want to identify the conserved quantity (energy) that can be associated with the system \eqref{eqn-for-psi-full}-\eqref{eqn-for-phi-full}. For that we apply Noether's Theorem which guarantees conservation of the extended energy given by
\begin{align}
\label{def-extended-energy}
\hat{E} &:= \langle \nabla_{\dt{\psi},\dt{\phi}} \mathcal{L}(\psi, \phi , \dt\psi , \dt\phi) , (\dt{\psi},\dt{\phi}) \rangle- \mathcal{L}(\psi,\phi,\dt{\psi},\dt{\phi}) \\
\nonumber&= \int_{\DOmega}  \frac{1}{2} |\nabla \phi|^2 +  \frac{\mu}{2} |\dt{\phi}|^2 + \frac{1}{2} (V[\rhophi] - \rhophi V^{\prime}[\rhophi]) \hspace{2pt}|\psi|^2 +  \frac{\mu \omega^2}{ 2} |\psi - \phi|^2 + \frac{1}{2} \rhophi V^{\prime}[\rhophi] \hspace{2pt}|2 \psi - \phi|^2  \hspace{2pt} dx.
\end{align}
Here we note that the value of the conserved (initial) energy depends on the initial condition for $\dt\phi(0)$, which we did not yet specify but which is required so that \eqref{eqn-for-phi-full} becomes a well-posed problem. There are at least two natural options for $\dt\phi(0)$, which we shall briefly discuss.

$\\$
{\bf Option 1.} One possibility is to use a compatibility condition to define $\dt\phi(0)$. Recall that we have $\psi(0)=\phi(0)=u_0$ and let $\rho_0:=|u_0|^2$. In this case, we observe by evaluating equation \eqref{eqn-for-psi-full} at time $t=0$ that
\begin{eqnarray*}
\ci \dt{\psi}(0) = - \frac{1}{2}\nabla^2 u_0 + \frac{1}{2} (V[\rho_0] - \rhophi V^{\prime}[\rho_0]) \psi  + \rho_0 V^{\prime}[\rho_0] u_0 =: \ci \dt{u}_0.
\end{eqnarray*}
This equation can be seen as a compatibility condition for smooth solutions. Hence, we could interpret $\dt{u}_0$ as a known quantity that can be defined (and computed) from $u_0$. Accordingly, we can exploit this knowledge and select the missing initial value as $\dt\phi(0)=\dt{u}_0$. However, from a computational point of view a different option can be more favorable.

$\\$
{\bf Option 2.} If the initial value for the first time derivative of $\phi$ is selected as $\dt{\phi}(0)=0$, then we observe that the preserved extended energy \eqref{def-extended-energy} is given by
\begin{align*}
\hat{E} = \frac{1}{2} \int_{\DOmega}  |\nabla u_0|^2 + V[\rho_0]  \hspace{2pt}|u_0|^2  \hspace{2pt} dx.
\end{align*}
In the case of the cubic nonlinearity $V[\rhophi] = 2 V_0 + \kappa \rhophi$ (which can be considered the most common model in the context of superfluidity), we see that the conserved value of $\hat{E}$ is identical to the energy associated with the original Lagrangian density, i.e., $E(u)=\int_{\DOmega} \frac{1}{2} |\nabla u| + V |u|^2 + \frac{\kappa}{2} |u|^4  \hspace{2pt} dx$. This motivates the choice $\dt{\phi}(0)=0$. In our numerical experiments we will show that this indeed yields very good approximations. 

As a side note, the issue of selecting a suitable initial value for $\dt{\phi}(0)$ is in fact closely related to finding slow manifolds of stiff mechanical systems, cf. \cite{AST12}.

$\\$
{\bf Preliminary notes on the adiabatic limit.} 
There is another interesting aspect about the conservation of $\hat{E}$ in the case $V[\rhophi] = 2 V_0 + \kappa \rhophi$. Assume that the trapping potential $V_0$ is nonnegative and that the particle interactions are repulsive, i.e., $\kappa>0$, then we can conclude that 
$$
\int_{\DOmega} \frac{\mu \omega^2}{ 2} |\psi - \phi|^2 \hspace{2pt} dx \le E(u_0).
$$
Hence, in the adiabatic limit $\mu \omega^2 \rightarrow \infty$, it must hold that $\int_{\DOmega}  |\psi - \phi|^2 \hspace{2pt} dx \rightarrow 0$. This shows a priori that $\psi$ and $\phi$ must converge in $L^2$ to the same function. Similarly, it can be also shown that the particle number is preserved for the limit. To see this, we multiply equation \eqref{eqn-for-psi-full} with $\psi^{\ast}$ and equation \eqref{eqn-for-phi-full} with $\phi^{\ast}/2$. Subtracting the arising terms and taking the imaginary part yields 
\begin{eqnarray*}
 \frac{1}{2} \frac{\mbox{d}}{\mbox{d}t} \int_{\DOmega} |\psi|^2 \hspace{2pt} dx &=&  \frac{\mu}{2} \frac{\mbox{d}}{\mbox{d}t}  \Im \int_{\DOmega} \dt{\phi} \phi^{\ast} \hspace{2pt} dx
\end{eqnarray*}
and therefore with $\dt{\phi}(0)=0$
\begin{eqnarray*}
 \int_{\DOmega} |\psi(t)|^2 \hspace{2pt} dx &=&  \int_{\DOmega} |\psi(0)|^2 \hspace{2pt} dx + \mu \Im \left( \int_{\DOmega}  \dt{\phi}(t) \phi^{\ast}(t) \hspace{2pt} dx \right).
\end{eqnarray*}
Together with the energy bound $\int_{\DOmega}\frac{1}{2}|\nabla \phi|^2 + \frac{\mu}{2} |\dt{\phi}|^2 \hspace{2pt} dx \le E(u_0)$ we conclude that
\begin{eqnarray*}
\left| \int_{\DOmega} |\psi(t)|^2 \hspace{2pt} dx-  \int_{\DOmega} |\psi(0)|^2 \hspace{2pt} dx  \right| \le  2  \mu \max_{s\le t} \| \dt{\phi}(s) \|_{L^2(\DOmega)} \| \phi^{\ast}(t) \|_{L^2(\DOmega)} \le C_{E(u_0)} \hspace{2pt} \mu^{1/2},
\end{eqnarray*}
for a constant $C_{E(u_0)}$ that only depends on $E(u_0)$ and $\DOmega$. In the adiabatic limit with $\mu \rightarrow 0$, we can hence conclude the conversation of mass as $\int_{\DOmega} |\psi(t)|^2 \hspace{2pt} dx \rightarrow  \int_{\DOmega} |\psi(0)|^2 \hspace{2pt} dx$.

\subsection{Adiabatic limit approximation}

In the next step, we want to approximate the adiabatic limit of the extended Euler-Lagrange equations \eqref{eqn-for-psi-full}-\eqref{eqn-for-phi-full}. For that we assume that the fictitious mass tends to zero, i.e., $\mu \rightarrow 0$, and that the frequency $\omega$ of the harmonic oscillator grows proportional to $\mu^{-1}$. In particular, we assume that $\mu \omega = {\rm constant}$ and hence $\omega \rightarrow \infty$. With this, the extended field variable $\phi$ is pushed closer and closer to the quantum state $\psi$. From the extended Born--Oppenheimer molecular dynamics it is known that the rate with which $\phi$ approaches $\psi$ is of order $|\psi-\phi| \sim {\cal O}(\omega^{-2})$ (cf. \cite{ANiklasson14}). In our numerical experiments in Section \ref{section:numerical-experiments} we observe the same relation for our extended Gross-Piatevskii model. Hence, we can formulate the adiabatic limit approximation (for $\mu, \omega^{-1} \rightarrow 0$) of the extended  Euler-Lagrange equations 
%\eqref{eqn-for-psi-full}-\eqref{eqn-for-phi-full}
\begin{eqnarray*}
\ci \dt{\psi} &=& - \frac{1}{2}\nabla^2 \psi + \frac{1}{2} (V[\rhophi] - \rhophi V^{\prime}[\rhophi]) \psi  + \rhophi V^{\prime}[\rhophi] (2 \psi - \phi) + \frac{\mu \omega^2}{2} ( \psi - \phi),\\
 \mu \ddt{\phi} &=&  - (V^{\prime}[\rhophi] + \rhophi V^{\prime\prime}[\rhophi]) |2 \psi - \phi|^2  \phi + \rhophi V^{\prime}[\rhophi] (2 \psi - \phi  ) + \rhophi  V^{\prime\prime}[\rhophi] \rho_{\psi} \phi + \mu \omega^2 (\psi - \phi)
\end{eqnarray*}
by asserting $\mu \rightarrow 0$; $\mu \omega = \mbox{constant}$ and $|\psi-\phi| \sim {\cal O}(\omega^{-2})$ to obtain the reduced model
\begin{eqnarray}
\label{final-sl-system} \ci \dt{\psi} &=& - \frac{1}{2}\nabla^2 \psi 
+ \frac{1}{2} V[\rhophi]  \psi  + \frac{1}{2} \rhophi V^{\prime}[\rhophi] (3 \psi - 2 \phi),\\
\nonumber \ddt{\phi} &=& \omega^2 (\psi - \phi).
\end{eqnarray}
As for the original extended system \eqref{eqn-for-psi-full}-\eqref{eqn-for-phi-full}, the equations are completed by the initial conditions
\begin{align*}
\psi(0)=\phi(0)=u_0 \qquad \mbox{and} \qquad \dt{\phi}(0)=0.
\end{align*}
For the rest of the paper, we shall assume that $\omega$ is large and we refer to the system \eqref{final-sl-system} as the shadow Lagrangian equations for superfluidity. Here we note that the system \eqref{final-sl-system} can be solved sequentially, where each discrete time step involves to first solve the wave equation $\ddt{\phi}= \omega^2 (\psi - \phi)$ explicitly and then the Schr\"odinger equation implicitly or semi-explicitly. 

As a final remark, we note that it is not recommended to reduce the term $(3 \psi -2 \phi)$ in \eqref{final-sl-system} to only $\psi$. Even though this could be justified with $|\psi-\phi| \sim {\cal O}(\omega^{-2})$, we made the numerical observation that the influence of $\phi$ in $ \rhophi V^{\prime}[\rhophi] (3 \psi - 2 \phi)$ has a significant stabilizing effect on corresponding discretizations.

\section{Shadow Lagrangian discretization}
\label{section:discretization}

In the next step, we will use the shadow Lagrangian formulation \eqref{final-sl-system} to introduce a new discretization for the classical Gross-Piatevskii equation. For that, we restrict our considerations from the previous section to the common case of cubic nonlinearities. More precisely, we let $V[\rho] = 2 V_0 + \kappa \rho$ for some potential $V_0 \in L^2(\DOmega)$ and some constant parameter $\kappa \in \mathbb{R}$. As mentioned before, the dynamics of Bose-Einstein condensates are described by the Gross-Pitaevskii equation where we seek the quantum state $u(\mathbf{x},t)$ of the condensate with
\begin{align}
\label{GPE}
\ci \hspace{2pt} \dt{u}(\mathbf{x},t) = -\frac{1}{2} \nabla^2 u(\mathbf{x},t) + V_0(\mathbf{x}) u(\mathbf{x},t) + \kappa |u(\mathbf{x},t)|^2 u(\mathbf{x},t) \quad \mbox{for } \mathbf{x} \in U, \enspace t \ge 0,
\end{align}
together with the initial condition $u(\mathbf{x},0)=u_0(\mathbf{x})$ and the boundary condition $u(\mathbf{x},t)=0$ for $\mathbf{x}\in\partial \DOmega$. Note that we denote the solution to \eqref{GPE} by $u$ to distinguish it from the quantum state $\psi$ that appears in the extended model. In order to find a numerical approximation to $u$, we first consider the shadow Lagrangian approximation $(\psi,\phi)$. In the Gross-Pitaevskii case, the shadow Lagrangian system \eqref{final-sl-system} simplifies to the equations
\begin{align}
\label{xGPE-1}
\ci \dt{\psi}(\mathbf{x},t)  &= -\frac{1}{2} \nabla^2 \psi (\mathbf{x},t) + V_0(\mathbf{x}) \psi(\mathbf{x},t) + \kappa |\phi(\mathbf{x},t)|^2 (2\psi(\mathbf{x},t) - \phi(\mathbf{x},t)) \quad \mbox{and}\\
\label{xGPE-2}
\ddt{\phi}(\mathbf{x},t)  &= \omega^2 (\psi(\mathbf{x},t) - \phi(\mathbf{x},t)),
\end{align}
for $\mathbf{x} \in U, \enspace t \ge 0$, 
together with the initial conditions $\psi(\mathbf{x},0)=\phi(\mathbf{x},0)=u_0(\mathbf{x})$ and $\dt{\phi}(\mathbf{x},0)=0$ and the boundary condition $\psi(\mathbf{x},t)=\phi(\mathbf{x},t)=0$ for $\mathbf{x}\in\partial \DOmega$. In order to discretize the shadow Lagrangian system in time, we introduce a small time step size $\tau>0$. As the frequency $\omega$ is supposed to grow to infinity in the adiabatic limit (and we wish to chose it as large as possible), we couple it to the time step size through the relation $\tau^{2} \omega^{2} = \beta_K$, or respectively 
$$
\omega = \sqrt{\beta_K} \tau^{-1}
$$
for some constant $\beta_K>0$ that we shall specify later. The coupling guarantees that $\omega\rightarrow \infty$ when we refine the time step size with $\tau \rightarrow 0$. Furthermore, the particular scaling $\tau \omega = \mbox{const}$ must be seen as a CFL (Courant-Friedrichs-Lewy) condition that ensures that an explicit leapfrog discretization of $\ddt{\phi} = \omega^2 (\psi- \phi)$ is stable. In fact, any coupling that would allow $\omega$ to grow faster, i.e., $\tau \omega^{1+\eps} = \mbox{const}$ for some $\eps>0$ is numerically unstable.

The first equation in the shadow Lagrangian formulation, i.e., equation \eqref{xGPE-1}, is discretized in time using a Crank-Nicolson-type approach. This choice is motivated by the observation that if we have consistency with $\psi=\phi$, then the discretization would be perfectly energy- and mass-conserving. The second equation,  i.e., equation \eqref{xGPE-2}, is discretized using a stabilized (dissipative) modification of the symplectic leapfrog (Verlet) scheme for the wave equation. This choice is motivated by the findings obtained in \cite{ANiklasson09,GZheng11} in the context of extended Lagrangian Born--Oppenheimer molecular dynamics. In conclusion we can formulate the following numerical method.
\begin{definition}[Dissipative Shadow Lagrangian Method]
\label{DSHM-formulation}
Let the initial values be given by $\psi^0=\phi^0 = u_0$ and let $K \in \mathbb{N}$ denote a dissipation order. Motivated by $\dt{\phi}(0)=0$, we also define the ghost values
$$
\phi^{-k} = \phi^0 \quad \mbox{for } k=1,\cdots,K+1.
$$
For simplicity, we denote averages by
\begin{align*}
\psi^{n+\frac{1}{2}} := \frac{\psi^{n+1} + \psi^n}{2} \qquad \mbox{and} \qquad
\phi^{n+\frac{1}{2}} := \frac{\phi^{n+1} + \phi^n}{2}.
\end{align*}
With this, the dissipative Shadow Lagrangian approximations $\psi^{n} \in H^1_0(U)$ and $\phi^{n}\in H^1_0(U)$ to $u$ at time $t^n$ are defined by
\begin{align}
\label{time-discretization-psi} \ci \psi^{n+1}  &= \ci \psi^n - \frac{\tau}{2} \nabla^2 \psi^{n+\frac{1}{2}} + \tau V \psi^{n+\frac{1}{2}} + \tau \kappa \frac{|\phi^{n+1}|^2 + |\phi^{n}|^2}{2} (2\psi^{n+\frac{1}{2}} - \phi^{n+\frac{1}{2}})
\end{align}
and
\begin{align}
\label{time-discretization-phi} \phi^{n+1}  &= 2\phi^{n} - \phi^{n-1} + \beta_K (\psi^n - \phi^n) + \alpha_K \sum_{k=0}^{K+1} c_k \phi^{n-k}.
\end{align}
The coefficients in \eqref{time-discretization-phi} depend on the selected order of dissipation $K$ and are exemplarily given according to the following table
\begin{center}
\begin{tabular}{|c||c|c|c|c|c|c|c|c|c|c|}
\hline $K$ & $\beta_K$ & $\alpha_K \times 10^{-3}$ & $c_0$ & $c_1$ & $c_2$ & $c_3$ & $c_4$ & $c_5$ & $c_6$ & $c_7$ \\
\hline
\hline 0  &   1.3   & 0     &        &  &  &  &   &  &  &   \\
\hline 2  &   1.69 & 150 &  -2   & 3   & 0 & -1 &   &  &  &   \\
\hline 3  &   1.75 & 57   &  -3   & 6   & -2 & -2 & 1  &  &  &   \\
\hline 4  &   1.82 & 18   &  -6   & 14 & -8 & -3 & 4 & -1 &  &   \\
\hline 5  &   1.84 & 5.5  &  -14 & 36 & -27 & -2 & 12 & -6 & 1 & \\
\hline 6  &   1.86 & 1.6  &  -36 & 99 & -88 & 11 & 32 & -25 & 8 & -1\\
\hline
\end{tabular}\end{center}
\end{definition}
There are some remarks that need to be made here. First, we recall that $\beta_K$ should be seen as the ratio that couples the time step size with the frequency, i.e., $\beta_K=\tau^2 \omega^{2}$. This means that \eqref{time-discretization-psi}-\eqref{time-discretization-phi} cannot be seen as a discretization of \eqref{xGPE-1}-\eqref{xGPE-2} for a fixed frequency, but rather a sequence of approximations for increasing frequencies. As such a sequence mimics the adiabatic limit process it is expected to converge to a solution of the original Gross-Pitaevskii model given by \eqref{GPE}. Hence, as desired, $\psi^{n}$ will be an approximation of the exact quantum state $u$ at time $t^n$.

Looking at the numbers in the table, we observe that for $K=0$, the discretization \eqref{time-discretization-phi} reduces to the standard symplectic leapfrog scheme without artificial dissipation. Even though we can expect this choice to work well for moderately large numbers of time steps, it is sensitive to an accumulation of numerical noise over time (entering through numerical errors in $\psi$). For $K=0$, this leads to severe deviations from the correct energy. We shall later confirm this numerically. In order to suppress the noise, a dissipative external force is added, which is modeled by the tail of the form $\alpha_K \sum_{k=0}^{K+1} c_k \phi^{n-k}$ in \eqref{time-discretization-phi}. This keeps the energy stable even over large numbers of time steps. As we will see in Section \ref{section:numerical-experiments}, any choice $K\ge3$ works very well in practice. Note however that the dissipation breaks time-reversibility and introduces damping, but only weakly in a higher odd order with respect to the integration time step size $\tau$. For $K \ge 4$, this would only be observed at very large time scales that are beyond the scope of the present study.
For a detailed derivation of the particular numbers for $\alpha_K$, $\beta_K$ and $c_k$ in the Born--Oppenheimer context, we refer to \cite{ANiklasson09}. Here we note that the optimization process that led to the values in the table is not a unique procedure and there exist many reasonable parameter combinations for each $K$ that all lead to reasonable discretizations.

Concerning equation \eqref{time-discretization-psi} it is worth to mention that if we replace $\phi^{n}$ by $\psi^{n}$ and $\phi^{n+1}$ by $\psi^{n+1}$ , then we obtain a fully energy conserving Crank-Nicolson discretization of \eqref{GPE} as it was studied and applied in \cite{ADK91,HeP17,San84}. Even though such a classical Crank-Nicolson discretization shows typically a high accuracy, it can be computationally demanding since it requires to solve a large nonlinear system of equations in every time step, which involves a repeated assembly of an updated stiffness matrix (within the iterations of a nonlinear solver). There is no such issue in our discretization \eqref{time-discretization-psi}, which only requires to solve a linear system in each time step.

Concerning the general computational complexity of the method, we note that the main cost per time step go (to roughly equal parts) to the assembly of the sparse system matrix that describes the system \eqref{time-discretization-psi} in a (standard) finite element space $V_h\subset H^1_0(U)$ and to the solving of the linear system afterwards. Here we note that the system matrix needs to be updated once per time step based on $\psi^n$, $\phi^n$ and $\phi^{n+1}$. Since the matrix is sparse this involves $\mathcal{O}(N)$ basic arithmetic operations, where $N$ is the dimension of $V_h$. The arising sparse system can be solved e.g. with an algebraic multigrid algorithm with a computational complexity of order $\mathcal{O}(\log(N) N)$. The cost for the explicit time stepping \eqref{time-discretization-phi} is computationally negligible. Finally we also note that the memory complexity per time step is of order $\mathcal{O}((K+1)N)$, where $K$ is the dissipation order.

As a final remark on the discretization we stress the important aspect that the difference between $\psi^n$ and $\phi^n$ can be seen as a consistency error. With respect to the $L^{\infty}(L^2)$-norm, the $L^{\infty}(H^1)$-norm and the error in energy difference it shows the same quadratic order convergence in $\tau$ as the error $\psi^n - u(t^n)$ itself. In particular, if the difference between $\psi^n$ and $\phi^n$ is large for a given $\tau$, then the discretization is not consistent and the error will be large as well. Hence, we can practically use the computable difference between $\psi^n$ and $\phi^n$ as an error indicator for the numerical method which allows us get an estimate for the size of the error between $\psi^n$ and $u(t^n)$. This is a significant practical advantage compared to alternative methods, where error indicators can be very hard to compute. In our numerical experiments in the next section, we will stress the usability of $\left| E(\psi^n) - E(\phi^n) \right|$ as an error indicator. On the one hand this indicator has the correct convergence order $\mathcal{O}(\tau^2)$ and on the other hand it  allows us to draw direct conclusions on the strength of the energy oscillations (which we want to keep as small as possible).

\section{Numerical experiments}
\label{section:numerical-experiments}

In our numerical experiments, we consider the Gross-Pitaevskii equation as given by \eqref{GPE}, where we fix the computational domain with $\DOmega=[-6,6]^2$ in all our experiments. The time discretization is according to the new Dissipative Shadow Lagrangian Method as stated in Definition \ref{DSHM-formulation}. The selected time step size $\tau$ is specified for each experiment individually. In order to study the influence of the dissipation order $K$ in \eqref{time-discretization-phi}, we will use the notation DS-$K$ to indicate which realization we used. For instance, DS-$K0$ will refer to \eqref{time-discretization-phi} with $K=0$ (i.e., the classical leapfrog discretization without dissipation), DS-$K2$ refers to \eqref{time-discretization-phi} with $K=2$, DS-$K3$ refers to \eqref{time-discretization-phi} with $K=3$; etc. The maximum number of time steps is denoted by $N$ and the maximum computing time by $T=t^N$ (we have $T=4$ for Model Problem 1 and $T=1$ for Model Problem 2). Consequently, $\psi^N$ denote the approximations at final time.

The space discretization in our experiments is based on piecewise linear Lagrange finite elements with a mesh width of approximately $h=0.05$. This space discretization is kept fixed in all our computations, including reference simulations. 

In order to quantify numerical errors, we also computed highly accurate reference solutions. These solutions were obtained with an energy- and mass- conservative Crank-Nicolson discretization of \eqref{GPE} that is known to be convergent to the exact quantum state $u$ \cite{HeP17}. Here we recall that the Crank-Nicolson discretization requires to solve nonlinear problems in each iteration, which makes it computationally more heavy. The time step size for the reference solution was chosen to be $\tau_{\mbox{\tiny ref}}=4\times 10^{-4}$ for Model Problem 1 in Section \ref{subsec-modprob-1} and  
%$\tau_{\mbox{\tiny ref}}=4\times 10^{-4}$ 
$\tau_{\mbox{\tiny ref}}=10^{-4}$ 
for Model Problem 2 in Section \ref{subsec-modprob-2}. For simplicity of the notation we will refer to this reference solution simply as $u$, since it is sufficiently close to the exact quantum state in \eqref{GPE}. Errors are measured in the $L^2$-norm and the $H^1$-norm which we recall as
$$
\| v\|_{L^2(\DOmega)} := \left(\int_{\DOmega} |v|^2 \hspace{2pt}dx\right)^{1/2} \quad \mbox{and} \quad 
\| v\|_{H^1(\DOmega)} := \left( \int_{\DOmega} |v|^2 + |\nabla v|^2  \hspace{2pt}dx\right)^{1/2} .
$$

\subsection{Model Problem 1}
\label{subsec-modprob-1}

%N64     -> 2^2/2^6 = 2^{-4}
%N128   -> 2^2/2^7 = 2^{-5}
%N256   -> 2^2/2^8 = 2^{-6}
%N512   -> 2^2/2^9 = 2^{-7}
%N1024 -> 2^2/2^10 = 2^{-8}

\begin{figure}[h!]
\centering
\includegraphics[scale=0.103]{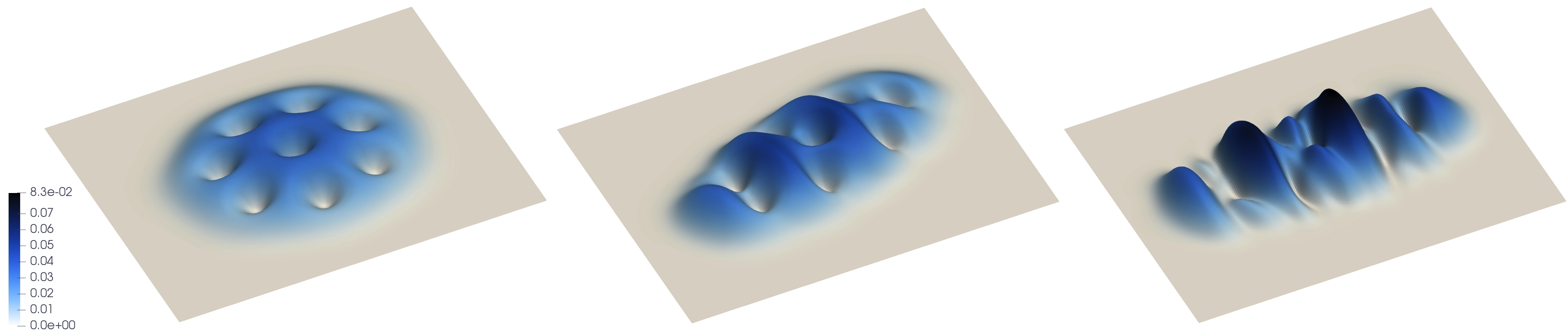}
\caption{\it Model problem 1. Left: ground state density $|u_0|^2$ at $t=0$. Middle: exact density $|u(t)|^2$ at time $t=2$ Right: exact density at final computing time $t=4$.}
\label{mod-prob-1-exact-u}
\end{figure}

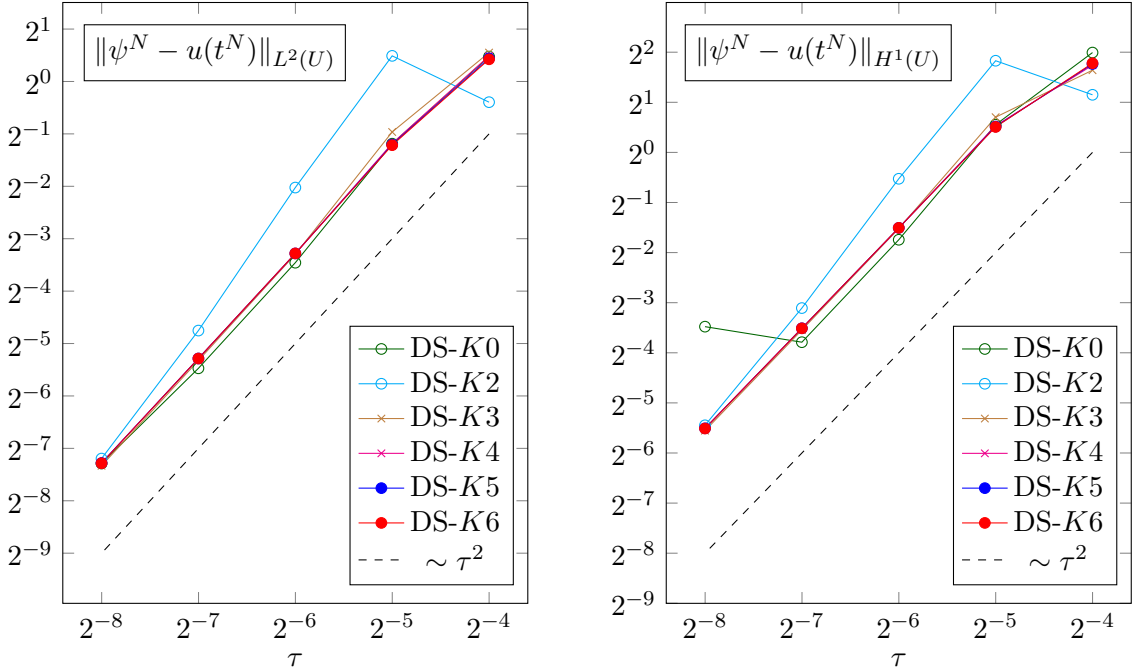
\begin{figure}
\begin{tikzpicture}
    \begin{axis}[
       legend pos=south east,%north west,
        height=0.62\textwidth,
        width=0.5\textwidth,
        %xmax   = 1,  xmin   = 0, ymax = 1, ymin = 0
        xlabel=$\tau$,%$\log_2(\tau)$,
        %ylabel=$\| \psi^N - u(t^N) \|_{L^2(\DOmega)}$,
        xmode=log, % use logarithmic scaling
        ymode=log,
        log basis y={2}, % basis of logarithm 
        log basis x={2} ] % basis of logarithm 
        %L2 errors
        %S0
        \addplot[mark=o,dark-green] coordinates {
        (2^-4, 1.39029) 
        (2^-5, 0.439258)
        (2^-6, 0.0910483)
        (2^-7, 0.0225274)
        (2^-8, 0.00639531)
        };
        % S3
        \addplot[mark=o,pantone_312] coordinates {
        (2^-4, 0.760627) 
        (2^-5, 1.40313)
        (2^-6, 0.245912)
        (2^-7, 0.0371147)
        (2^-8, 0.00683101)
        };    
        % S4
        \addplot[mark=x, brown] coordinates {
        (2^-4, 1.46557) 
        (2^-5, 0.512817)
        (2^-6, 0.101773)
        (2^-7, 0.0249763)
        (2^-8, 0.00621105)
        };
         % S5
        \addplot[mark=x,magenta] coordinates {
        (2^-4, 1.40169) 
        (2^-5, 0.440617)
        (2^-6, 0.102592)
        (2^-7, 0.0256371)
        (2^-8, 0.00643092)
        };
        % S6
        \addplot[mark=*, blue] coordinates {
        (2^-4, 1.35811) 
        (2^-5, 0.434221)
        (2^-6, 0.10291)
        (2^-7, 0.0256608)
        (2^-8, 0.00644343)
        }; 
        % S7
        \addplot[mark=*, red] coordinates {   
        (2^-4, 1.34051) 
        (2^-5, 0.430766)
        (2^-6, 0.102929)
        (2^-7, 0.0256422)
        (2^-8, 0.00642247)
        };
        % \tau^2
        \addplot[black, dashed] coordinates {
        (2^-4,2^-1) 
        (2^-5,2^-3) 
        (2^-6,2^-5) 
        (2^-7,2^-7)
        (2^-8,2^-9) 
        }; 
        \legend{%
DS-$K0$,
DS-$K2$,
DS-$K3$,
DS-$K4$,
DS-$K5$,
DS-$K6$,
$\sim \tau^{2}$,
%Component \textbf{d},
}
    \end{axis}
   \node[draw] at (2,7.3) {$\| \psi^N - u(t^N) \|_{L^2(\DOmega)}$};
\end{tikzpicture}
\hspace{20pt}
\begin{tikzpicture}
    \begin{axis}[
       legend pos=south east,%north west,
        height=0.62\textwidth,
        width=0.5\textwidth,
        %xmax   = 1,  xmin   = 0, ymax = 1, ymin = 0
        xlabel=$\tau$,%$\log_2(\tau)$,
        %ylabel=$\| \psi^N - u(t^N)\|_{H^1(\DOmega)}^{\mbox{\rm\tiny rel}}$,
        xmode=log, % use logarithmic scaling
        ymode=log,
        log basis y={2}, % basis of logarithm 
        log basis x={2} ] % basis of logarithm 
        %H1 errors
        %S0
        \addplot[mark=o,dark-green] coordinates {
        (2^-4, 1.601 * 2.48542) 
        (2^-5, 0.590507 * 2.48542) 
        (2^-6, 0.120171 * 2.48542) 
        (2^-7, 0.0291315 *2.48542) 
        (2^-8, 0.0361147*2.48542)  
        };
        % S3
        \addplot[mark=o,pantone_312] coordinates {
        (2^-4, 0.892576*2.48542) 
        (2^-5, 1.42981*2.48542) 
        (2^-6, 0.279347*2.48542) 
        (2^-7, 0.0467019*2.48542) 
        (2^-8, 0.00922885 * 2.48542) 
        };    
        % S4
        \addplot[mark=x, brown] coordinates {
        (2^-4, 1.25432 * 2.48542) 
        (2^-5, 0.6548280 * 2.48542) 
        (2^-6, 0.140362 * 2.48542) 
        (2^-7, 0.0344798 * 2.48542) 
        (2^-8, 0.00856552 * 2.48542) 
        };
         % S5
        \addplot[mark=x,magenta] coordinates {
        (2^-4, 1.344 * 2.48542) 
        (2^-5, 0.582326 * 2.48542) 
        (2^-6, 0.141309 * 2.48542) 
        (2^-7, 0.0352943 * 2.48542) 
        (2^-8, 0.00883575 * 2.48542) 
        };
        % S6
        \addplot[mark=*, blue] coordinates {
        (2^-4, 1.36751 * 2.48542) 
        (2^-5, 0.575602 * 2.48542) 
        (2^-6, 0.141674 * 2.48542) 
        (2^-7, 0.0353216 * 2.48542) 
        (2^-8, 0.00885083 * 2.48542) 
        }; 
        % S7
        \addplot[mark=*, red] coordinates {   
        (2^-4, 1.38122 * 2.48542) 
        (2^-5, 0.572128 * 2.48542) 
        (2^-6, 0.141686 * 2.48542) 
        (2^-7, 0.0352915 * 2.48542) 
        (2^-8, 0.00882518 * 2.48542) 
        };
        % \tau^2
        \addplot[black, dashed] coordinates {
        (2^-4,2^0) 
        (2^-5,2^-2) 
        (2^-6,2^-4) 
        (2^-7,2^-6)
        (2^-8,2^-8) 
        }; 
        \legend{%
DS-$K0$,
DS-$K2$,
DS-$K3$,
DS-$K4$,
DS-$K5$,
DS-$K6$,
$\sim \tau^{2}$,
%Component \textbf{d},
}
    \end{axis}
    \node[draw] at (2,7.3) {$\| \psi^N - u(t^N)\|_{H^1(\DOmega)}$}; %^{\mbox{\rm\tiny rel}}
\end{tikzpicture}
\caption{{\it Model Problem 1. Comparison of accuracies for the Dissipative Shadow Lagrangian Method with different dissipation orders $K$. Notably, the accuracies for $K=3,4,5,6$ are extremely close to each other.}}
\label{model-problem-1-dissipation-error}
\end{figure}

In the first numerical experiment we consider a setup that involves a rotating Bose-Einstein condensate (BEC) with repulsive particle interactions characterized by the parameter $\kappa=100$. In the following, $V_{\mbox{\tiny har}}^{\boldsymbol{\gamma}}$ denotes a harmonic trapping potential of the form
$$V_{\mbox{\tiny har}}^{\boldsymbol{\gamma}}(\mathbf{x}) = \frac{1}{2} ( \gamma_1 x_1^2 + \gamma_2 x_2^2),$$
where $\gamma_1>0$ and $\gamma_2>0$ denote the trapping frequencies. A stirring potential is modeled by the angular momentum operator $\mathcal{L}_z = - \ci \left( x \partial_y - y \partial_x \right)$ (i.e., we have a rotation around the $z$-axis). Fixing the angular velocity with $\Omega=0.8$ and the trapping frequencies with $\gamma_1=\gamma_2=1$ we start with computing a corresponding vortex ground state $u_0$ which is an $L^2$-normalized minimizer of the energy 
$$
E_0(v) = \int_{\DOmega}  \frac{1}{2}|\nabla v|^2 + V_{\mbox{\tiny har}}^{\boldsymbol{\gamma}}|v|^2 + \Omega \mathcal{L}_z(v) v^{\ast} + \frac{\kappa}{2} |v|^4 \hspace{2pt}dx, \qquad
(\mbox{where } \gamma_1 = \gamma_2  =1)
$$
among all states $v \in H^1_0(\DOmega)$ with $\| v\|_{L^2(\DOmega)}=1$. Ground states can be computed by solving a nonlinear eigenvalue problem (cf. \cite{AHP19,BaD04,BWM05,DaK10,HeP18,JKM14}). In our experiment, we computed a ground state density that contains eight vortices with a quantized circulation. It is depicted in Figure \ref{mod-prob-1-exact-u} (left). The vortices are density singularities and arise commonly in rotating BEC's due to the superfluid character of the condensate (cf. \cite{FSF01,MCB01,MCW00}). The ground state $u_0$ is used as a starting value for the time-dependent problem \eqref{GPE}. In order to trigger some interesting dynamics, the stirring potential is switched off and the harmonic potential is reconfigured to an anisotropic trap with new trapping frequencies $\gamma_1=2$ and $\gamma_2=1$. We observe that the condensate continues to rotate, but the initial rotational symmetry is broken and it deforms quickly into cigar-shaped object. The dynamics of the density $|u(t)|^2$ at snapshot times $t=0$, $t=2$ and final time $T=4$ are shown in Figure \ref{mod-prob-1-exact-u}. Important time-invariants are the conservation of mass $m(u):=\| u \|_{L^2(\DOmega)}=1$ and energy
\begin{align}
\label{energy-mod-prob-1}
E(u) := \int_{\DOmega}  \frac{1}{2}|\nabla u|^2 + V_{\mbox{\tiny har}}^{\boldsymbol{\gamma}} |u|^2 + \frac{\kappa}{2} |u|^4 \hspace{2pt}dx,
\quad
\mbox{where } \gamma_1 = 2 \mbox{ and } \gamma_2 =1.
\end{align}
In the first experiment, we investigate the influence of the dissipation order $K$ in the Dissipative Shadow Lagrangian Method (DS-$K$ as given by Definition \ref{DSHM-formulation}). The corresponding results are depicted in Figure \ref{model-problem-1-dissipation-error}, both for the $L^2$-error and the $H^1$-error at final computing time $t^N=4$. First, we note that all realizations converge as expected with a quadratic order in the step size, i.e., with order $\mathcal{O}(\tau^2)$. Next, we observe that the simple discretization without artificial dissipation, i.e., DS-$K0$, performs well for large step sizes $\tau$. However, for $\tau=2^{-8}$ an accumulation of numerical errors becomes visible, which puts this choice behind the realizations with artificial dissipative force. Even though not depicted in the graph of Figure \ref{model-problem-1-dissipation-error}, we observed that the poor performance of the DS-$K0$ becomes even more pronounced for $\tau=2^{-9}$ and the method becomes unstable as long as we keep the ratio $\tau^2 = 1.3 \times \omega^{-2}$ fixed. Among the dissipative methods, 
the realization with lowest dissipation order, i.e., DS-$K2$, shows the least accuracy. However, for the higher orders  $K=3,4,5,6$, there are no longer any significant differences and any choice leads to a method with a very good performance and no signs of instabilities. This was part of a general pattern and we made the same observations in other experiments. Hence, for the rest of the paper, we will only consider DS-$K5$ as representative scheme for the general accuracy (for $K\ge3$).

% Model Problem 1 - Comparison dissipation orders

%\input{model-problem-1-dissipation-error}

%N64     -> 2^2/2^6 = 2^{-4}
%N128   -> 2^2/2^7 = 2^{-5}
%N256   -> 2^2/2^8 = 2^{-6}
%N512   -> 2^2/2^9 = 2^{-7}
%N1024 -> 2^2/2^10 = 2^{-8}

% Model Problem 1 - 

%\input{model-problem-1-comparison-and-consistency}

%N64     -> 2^2/2^6 = 2^{-4}
%N128   -> 2^2/2^7 = 2^{-5}
%N256   -> 2^2/2^8 = 2^{-6}
%N512   -> 2^2/2^9 = 2^{-7}
%N1024 -> 2^2/2^10 = 2^{-8}

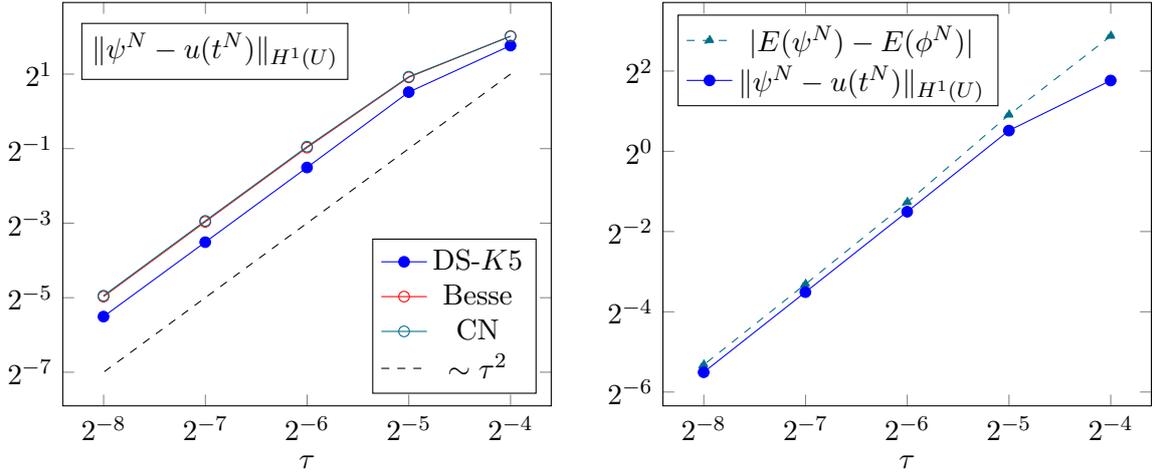
\begin{figure}
% using S6 absolute errors
\begin{tikzpicture}
    \begin{axis}[
       legend pos=south east,%north west,
        height=0.45\textwidth,
        width=0.52\textwidth,
        %xmax   = 1,  xmin   = 0, ymax = 1, ymin = 0
        xlabel=$\tau$,%$\log_2(\tau)$,
        %ylabel=$\| \psi^N- u ( t^N)\|_{H^1(\DOmega)}$,
        xmode=log, % use logarithmic scaling
        ymode=log,
        log basis y={2}, % basis of logarithm 
        log basis x={2} ] % basis of logarithm 
        \addplot[mark=*,blue] coordinates {
        % S6 = SH-K5
        (2^-4,3.39884) % N64
        (2^-5,1.43061) %N128
        (2^-6,0.35212) %N256
        (2^-7,0.0877892) %N512
        (2^-8,0.0219981)}; %N1024
        %
        % Besse:
        \addplot[mark=o,red] coordinates {
        (2^-4,4.05285) % N64
        (2^-5,1.88319) %N128
        (2^-6,0.509708) %N256
        (2^-7,0.128183) %N512
        (2^-8,0.0319107)}; %N1024
        %
         % CN:
        \addplot[mark=o, pantone_315] coordinates {
        (2^-4,4.03707) % N64
        (2^-5,1.90593) %N128
        (2^-6,0.518699) %N256
        (2^-7,0.130504) %N512
        (2^-8,0.0324759)}; %N1024
        \addplot[black, dashed] coordinates {
        (2^-4,2^1) 
        (2^-5,2^-1) 
        (2^-6,2^-3) 
        (2^-7,2^-5)
        (2^-8,2^-7)};
        \legend{%
DS-$K5$,
Besse,
CN,
$\sim \tau^{2}$,
}
    \end{axis}
        \node[draw] at (2,4.7) {$\| \psi^N- u ( t^N)\|_{H^1(\DOmega)}$};
\end{tikzpicture}
\hspace{10pt}
% S6
\begin{tikzpicture}
    \begin{axis}[
       legend pos=north west,
        height=0.45\textwidth,
        width=0.52\textwidth,
        %xmax   = 1,  xmin   = 0, ymax = 1, ymin = 0
        xlabel=$\tau$,%$\log_2(\tau)$,
        %ylabel=$\log_2(\| \psi - \phi \|)$,
        xmode=log, % use logarithmic scaling
        ymode=log,
        log basis y={2}, % basis of logarithm 
        log basis x={2} ] % basis of logarithm 
        \addplot[mark=triangle*,dashed,pantone_315] coordinates {
        (2^-4,7.31488) % N64
        (2^-5,1.87504) %N128
        (2^-6,0.41398) %N256
        (2^-7,0.100666) %N512
        (2^-8,0.0249971)}; %N1024
        \addplot[mark=*,blue] coordinates {
        (2^-4,3.39884) % N64
        (2^-5,1.43061) %N32
        (2^-6,0.35212) %N64
        (2^-7,0.0877892) %N128
        (2^-8,0.0219981)};  %N256
        \legend{%
$|E(\psi^N) - E(\phi^N)|$,%$ (t,\cdot)- \phi(t,\cdot) \|_{L^2(\DOmega)}$,
$\| \psi^N- u ( t^N)\|_{H^1(\DOmega)}$,
%Component \textbf{d},
}
    \end{axis}
\end{tikzpicture}
\caption{{\it Model Problem 1. Left: Comparison of the $H^1$-errors for DS-$K5$, the classical Crank-Nicolson scheme and the Besse-relaxation scheme. Right: comparison of the exact $H^1$-error for DS-$K5$ with the estimated error using the consistency indicator $|E(\psi^N) - E(\phi^N)|$.}}
\label{model-problem-1-comparison-and-consistency}
\end{figure}

Next, we study how the accuracy of the new approach compares to other established schemes. In this comparison we consider the Crank-Nicolson method (CN) which is obtained from \eqref{time-discretization-psi} if we replace $\phi^{n}$ and $\phi^{n+1}$ by $\psi^n$ and $\psi^{n+1}$. This method is known to conserve the exact mass $m(u)$ and energy $E(u)$ as given by \eqref{energy-mod-prob-1}, but it requires to solve a nonlinear problem in each time step. Hence, it is more expensive than DS-$K$. The second approach in our comparison is the Besse relaxation scheme proposed in \cite{Bes04}, which also conserves energy and mass and which only requires to solve linear problems in each step. This method yields a good reference as it is known to be extremely efficient and accurate (cf. the detailed numerical studies in \cite{HeW19}). Its computational complexity is the same as the complexity of the newly proposed DS-$K$ approach. The results are depicted in Figure \ref{model-problem-1-comparison-and-consistency} (left), exemplarily for the $H^1$-error as the picture for the $L^2$-error was the same. We can see that all methods show a quadratic convergence of order $\mathcal{O}(\tau^2)$. The errors for the CN and Besse relaxation method are almost identical, where the error for DS-$K5$ is around $30\%$ smaller. This shows that our new approach is competitive with existing approaches.

\input{model-problem-1-energy-mass-plots}

Let us next recall that the difference between $\psi^N$ and the auxiliary wave function $\phi^N$ can be seen as a consistency error that converges with order $\omega^{-2}= \beta_K^{-1} \tau^2$ to zero, i.e., with the same order as the error $\| \psi^N - u(t^N) \|_{H^1(\DOmega)}$ itself. Motivated by this observation, we suggest to use the energy consistency $\eta^N:=|E(\psi^N) - E(\phi^N)|\sim \tau^2$ as an error indicator. Due to the initial consistency $\psi(0)=\phi(0)=u_0$, the error indicator $\eta^N$ measures how well $\phi^N$ fits the equation for $\psi^N$ (note that $\psi^N=\phi^N=u(T)$ in the adiabatic limit $\omega\rightarrow \infty$) and the strength of the arising energy deviation. In Figure \ref{model-problem-1-comparison-and-consistency} (right) we plot the consistency error $\eta^N$ against the exact error  $\| \psi^N - u(t^N) \|_{H^1(\DOmega)}$ for DS-$K5$. We observe that both converge with the order $\mathcal{O}(\tau^2)$ and that $\eta^N$ is very close and always slightly above the real error. Hence, $\eta^N$ resembles nicely the overall accuracy.

Finally, we study the evolution of mass $m$ and energy $E$ (see \eqref{energy-mod-prob-1}) for the DS-$K$ approximations. In the original Gross-Pitaevskii model \eqref{GPE}, both quantities are conserved over time. However, this is no longer the case when switching to the extended shadow Lagrangian model \eqref{xGPE-1}-\eqref{xGPE-2}. Therefore it is important to ensure that the deviation from both quantities does not become too large and that the amplitude of the variations reduces with $\tau \to 0$. The corresponding results are shown in Figure \ref{model-problem-1-energy-mass-plots}, for DS-$K5$ with $\tau=2^{-4},2^{-5},2^{-6}$. We observe that both mass and energy are oscillating, where the strength of the oscillations damps out with time. The oscillations are strongest for $\tau=2^{-4}$, but are barely visible for $\tau=2^{-6}$ where we are already close to an ideal conservation of both quantities. In fact, the nature of the variations is very similar to what is typically observed for symplectic approximations of \eqref{GPE} (cf. \cite{HeW19}). As a very interesting observation, we see that the graphs of the energy and mass evolution follow each other very closely. This indicates that any loss in mass correlates with a simultaneous loss of energy.

\subsection{Model Problem 2}
\label{subsec-modprob-2}

\begin{figure}[h!]
\centering
\includegraphics[scale=0.14]{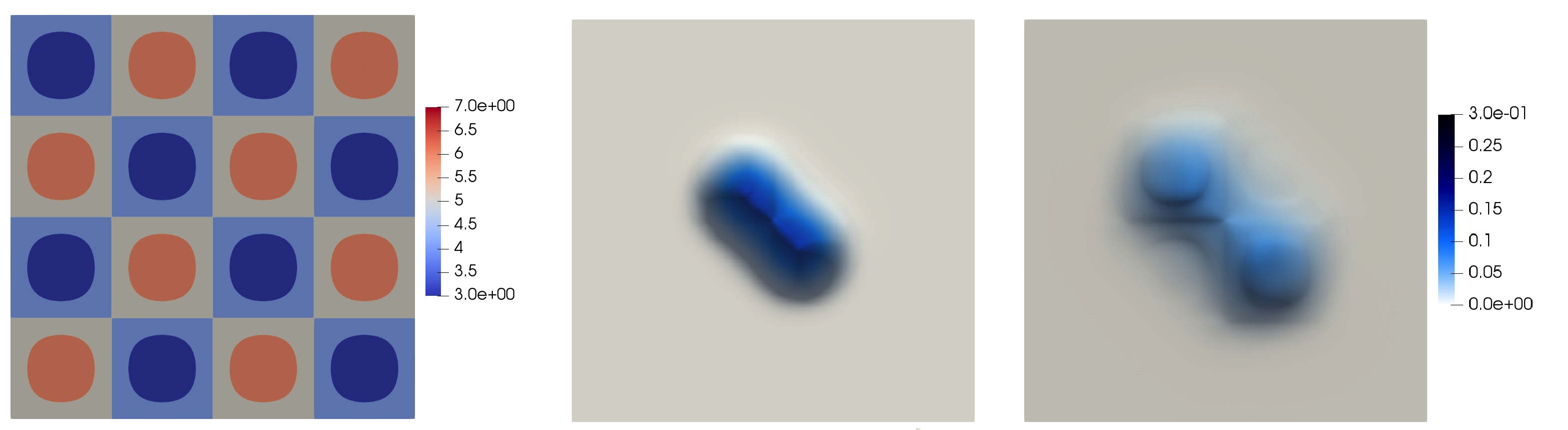}
\caption{\it Model problem 2. Left: potential $V_0$ as given by \eqref{V0-modprob2}. Middle: ground state density $|u_0|^2$, i.e., $|u(t)|^2$ at $t=0$. Right: exact density $|u(t)|^2$ at final computing time $t=1$.}
\label{mod-prob-2-pot-and-exact-u}
\end{figure}

In the second model problem we consider the setting that the potential is given by a discontinuous checker-board potential of the form
\begin{align}
\label{V0-modprob2}
V_0(\mathbf{x})= \left\lfloor5 + 2 \sin( \frac{\pi}{3} x_1 ) \sin( \frac{\pi}{3} x_2 ) \right\rfloor,
\end{align}
where $\lfloor\cdot \rfloor$ denotes the floor function. The potential is depicted in Figure \ref{mod-prob-2-pot-and-exact-u} (left). With this choice we want to trigger a loss of regularity of the quantum state $u$ which in consequence causes a drop in the expected convergence rates (cf. \cite{HeP17}). The particle interactions are repulsive where we select the corresponding parameter to be $\kappa=20$. It only remains to specify the initial value $u_0$ in \eqref{GPE}. Here we select $u_0$ to be the non-negative and $L^2$-normalized ground state associated with the initial energy
$$
E_0(v) = \int_{\DOmega}  \frac{1}{2}|\nabla v|^2 + (V_0 + V_{\mbox{\tiny har}})|v|^2 + \frac{\kappa_0}{2} |v|^4 \hspace{2pt}dx,
$$
where $\kappa_0=10$ and $V_{\mbox{\tiny har}}$ is the harmonic trapping potential with $V_{\mbox{\tiny har}}(\mathbf{x}) = \frac{1}{2} (x_1^2 + x_2^2)$. The corresponding ground state density $|u_0|^2$ is visualized in Figure \ref{mod-prob-2-pot-and-exact-u} (middle). We emphasize that the harmonic potential $V_{\mbox{\tiny har}}$ is switched off when studying the dynamics and that $\kappa$ changes from $10$ to $20$. Hence, the conserved quantity of $u$ is the energy
\begin{align}
\label{energy-mod-prob-2}
E(u) = \int_{\DOmega}  \frac{1}{2}|\nabla u|^2 + V_0 |u|^2 + \frac{\kappa}{2} |u|^4 \hspace{2pt}dx.
\end{align}

The numerical comparisons for Model Problem 2 are carried out at final time $t^N=1.0$. The exact density of the quantum state $u$ at $t^N$ is depicted in Figure \ref{mod-prob-2-pot-and-exact-u} (right). We observe that the density slowly evolves into the shape of the discontinuous trapping potential $V_0$.

% Model Problem 2 - 

%\input{model-problem-2-comparison}

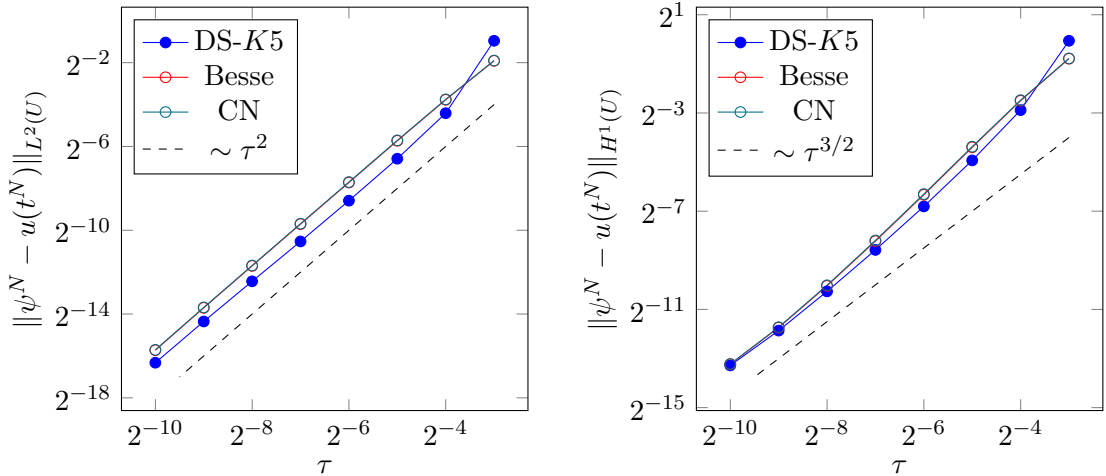
\begin{figure}
% using S6 absolute errors
\begin{tikzpicture}
    \begin{axis}[
       legend pos=north west,
        height=0.45\textwidth,
        width=0.45\textwidth,
        %xmax   = 1,  xmin   = 0, ymax = 1, ymin = 0
        xlabel=$\tau$,%$\log_2(\tau)$,
        ylabel=$\| \psi^N- u ( t^N)\|_{L^2(\DOmega)}$,
        xmode=log, % use logarithmic scaling
        ymode=log,
        log basis y={2}, % basis of logarithm 
        log basis x={2} ] % basis of logarithm 
        \addplot[mark=*,blue] coordinates {
        (2^-3,0.515699)
        (2^-4,0.0469697) 
        (2^-5,0.0104319) 
        (2^-6,0.00260111) 
        (2^-7,0.000674903)
        (2^-8,0.000180505) 
        (2^-9,4.77786e-05) 
        (2^-10,1.22116e-05)};
        % Besse
        \addplot[mark=o,red] coordinates {
        (2^-3,0.265553) 
        (2^-4,0.073295) 
        (2^-5,0.0188561) 
        (2^-6,0.00474125) 
        (2^-7,0.00119648)
        (2^-8,0.00030029) 
        (2^-9,7.5064e-05) 
        (2^-10,1.84522e-05)
        };    
         % CN
        \addplot[mark=o, pantone_315] coordinates {
        (2^-3,0.268098) 
        (2^-4,0.0742525) 
        (2^-5,0.0191304) 
        (2^-6,0.00481266) 
        (2^-7,0.00121453) 
        (2^-8,0.000304798) 
        (2^-9,7.61627e-05)
        (2^-10,1.86986e-05)
        };    
        %RK
        %\addplot[mark=o, pantone_312] coordinates {   
        % (2^-3,0.396235) 
        % (2^-4,0.113833) 
         % (2^-5,0.0296107) 
         % (2^-6,0.00746798) 
         % (2^-7,0.00187905) 
         % (2^-8,0.000470524) 
         % (2^-9,0.000117417) 
         % (2^-10,2.89286e-05)
         %};    
        \addplot[black, dashed] coordinates {
        (2^-3,2^-4) 
        (2^-4,2^-6) 
        (2^-5,2^-8) 
        (2^-6,2^-10) 
        (2^-7,2^-12)
        (2^-8,2^-14) 
        (2^-9,2^-16) 
        (2^-9.5,2^-17)}; 
        \legend{%
DS-$K5$,
Besse,
CN,
%RK,
$\sim \tau^{2}$,
%Component \textbf{d},
}
    \end{axis}
\end{tikzpicture}
\hspace{10pt}
%
% absolute errors
\begin{tikzpicture}
    \begin{axis}[
       legend pos=north west,
        height=0.45\textwidth,
        width=0.45\textwidth,
        %xmax   = 1,  xmin   = 0, ymax = 1, ymin = 0
        xlabel=$\tau$,%$\log_2(\tau)$,
        ylabel=$\| \psi^N- u ( t^N)\|_{H^1(\DOmega)}$,
        xmode=log, % use logarithmic scaling
        ymode=log,
        log basis y={2}, % basis of logarithm 
        log basis x={2} ] % basis of logarithm 
        \addplot[mark=*,blue] coordinates {
        (2^-3,0.958906) %N8
        (2^-4,0.13544) %N16
        (2^-5,0.0327752) %N32
        (2^-6,0.00894328) %N64
        (2^-7,0.0026185) %N128
        (2^-8,0.000812523)  %N256
        (2^-9,0.000268605) %N512 %0.000270296
        (2^-10,0.000100553)}; %1024
        % Besse
        \addplot[mark=o,red] coordinates {
        (2^-3,0.578454)
        (2^-4,0.176729) 
        (2^-5,0.0473946) 
        (2^-6,0.0124156) %N64 
        (2^-7,0.00335802)
        (2^-8,0.000950263) 
        (2^-9,0.00029238) 
        (2^-10,0.000104046)
        };    
         % CN
        \addplot[mark=o, pantone_315] coordinates {
        (2^-3,0.583441)
        (2^-4,0.179429) 
        (2^-5,0.0482962) 
        (2^-6,0.0126641) %N64
        (2^-7,0.00342158)
        (2^-8,0.000965961)
        (2^-9,0.000295851) 
        (2^-10,0.000104717)
        };    
        \addplot[black, dashed] coordinates {
        (2^-3,2^-4) %
        (2^-4,2^-5.5)  % 
        (2^-5,2^-7) 
        (2^-6,2^-8.5) 
        (2^-7,2^-10)
        (2^-8,2^-11.5) 
        (2^-9,2^-13)
        (2^-9.5,2^-13.75)}; 
        %
        % (2^-3,2^-4) %
        % (2^-4,2^-6)  % 
        % (2^-5,2^-8) 
        % (2^-6,2^-10) 
        % (2^-7,2^-12)
        % (2^-8,2^-14) 
        % (2^-9,2^-16) 
        % (2^-9.5,2^-17)}; 
        \legend{%
DS-$K5$,
Besse,
CN,
%$\sim \tau^{2}$,
$\sim \tau^{3/2}$,
}
    \end{axis}
\end{tikzpicture}
\caption{{\it Model Problem 2. Comparison of the accuracies for DS-$K5$, the classical Crank-Nicolson scheme and the Besse-relaxation scheme. Left: comparison of $L^2$-errors. Right: comparison of $H^1$-errors with visibly reduced convergence rates.}}
\label{model-problem-2-comparison}
\end{figure}

We start with a comparison of the numerical errors for DS-K$5$, the Crank-Nicolson scheme and the popular Besse relaxation scheme (cf. Section \ref{subsec-modprob-1} for references). The results are depicted in Figure \ref{model-problem-2-comparison}. For the $L^2$-error we observe an optimal convergence rate of order $\mathcal{O}(\tau^{2})$ for all three methods. As in Model Problem 1, the error graphs for the Crank-Nicolson and Besse method are basically indistinguishable, whereas the errors for DS-K$5$ are slightly smaller, reconfirming the competitive performance. More interestingly are the obtained $H^1$-errors which suffer from degenerate  convergence rates due to the loss of regularity triggered by the discontinuous potential. In the right graph of Figure \ref{model-problem-2-comparison} we can see that the convergence speed is consistently below $\mathcal{O}(\tau^{2})$ and reduces asymptotically to around $\mathcal{O}(\tau^{3/2})$ for all three methods.

% Model Problem 2 - 

%\input{model-problem-2-consistency-errors}

\begin{figure}
% convergence consistency error
\begin{tikzpicture}
    \begin{axis}[
       legend pos=north west,
        %height=0.6\textwidth,
        %width=0.6\textwidth,
        %xmax   = 1,  xmin   = 0, ymax = 1, ymin = 0
        xlabel=$\tau$,%$\log_2(\tau)$,
        %ylabel=$\log_2(\| \psi - \phi \|)$,
        xmode=log, % use logarithmic scaling
        ymode=log,
        log basis y={2}, % basis of logarithm 
        log basis x={2} ] % basis of logarithm 
        \addplot[mark=triangle*,dashed,red] coordinates {
        (2^-4,0.155131) 
        (2^-5,0.045649) 
        (2^-6,0.00524164) 
        (2^-7,0.00121233)
        (2^-8,0.000303237) 
        (2^-9,7.57949e-05) 
        (2^-10,1.8947e-05)};
        \addplot[mark=triangle*,dashed,blue] coordinates {
        (2^-4,0.275113) %absolute: 0.275113 relative: 0.168489
        (2^-5,0.068801 ) %absolute: 0.068801 relative: 0.041412457
        (2^-6,0.0103981)  %absolute: 0.0103981 relative: 0.0062491
        (2^-7,0.0025794) %absolute: 0.0025794 relative: 0.0015501
        (2^-8, 0.000650533)  %absolute: 0.000650533 relative: 0.00039096
        (2^-9,0.000165937)  %absolute: 0.000165937 relative: 0.0000997259
        (2^-10,4.20086e-05)}; %absolute: 4.20086e-05 relative: 0.0000252466
        \addplot[black, dashed] coordinates {
        (2^-4,2^-6) 
        (2^-5,2^-8) 
        (2^-6,2^-10) 
        (2^-7,2^-12)
        (2^-8,2^-14) 
        (2^-9,2^-16)}; 
        \legend{%
$\| \psi^N - \phi^N \|_{L^2(\DOmega)}$,%$ (t,\cdot)- \phi(t,\cdot) \|_{L^2(\DOmega)}$,
$\| \psi^N - \phi^N \|_{H^1(\DOmega)}$,%$\| \psi(t,\cdot) - \phi(t,\cdot) \|_{H^1(\DOmega)}$,
$\sim \tau^{2}$,
%Component \textbf{d},
}
    \end{axis}
\end{tikzpicture}
\begin{tikzpicture}
    \begin{axis}[
       legend pos=north west,
        %height=0.6\textwidth,
        %width=0.6\textwidth,
        %xmax   = 1,  xmin   = 0, ymax = 1, ymin = 0
        xlabel=$\tau$,%$\log_2(\tau)$,
        %ylabel=$\log_2(\| \psi - \phi \|)$,
        xmode=log, % use logarithmic scaling
        ymode=log,
        log basis y={2}, % basis of logarithm 
        log basis x={2} ] % basis of logarithm 
        \addplot[mark=triangle*,dashed,pantone_315] coordinates {
        (2^-4, 0.887886) %N16
        (2^-5, 0.298352) %N32
        (2^-6, 0.0585458) %N64
        (2^-7, 0.0135952) %128
        (2^-8, 0.00342122) %256
        (2^-9, 0.000855114) %512
        (2^-10, 0.000213765)}; %1024
        \addplot[mark=*,blue] coordinates {
        (2^-4,0.13544) %N16
        (2^-5,0.0327752) %N32
        (2^-6,0.00894328) %N64 %EOC: 1.8737
        (2^-7,0.0026185) %N128 %EOC: 1.772
        (2^-8,0.000812523)  %N256  %EOC: 1.6882
        (2^-9,0.000268605) %N512 %EOC: 1.5969
        (2^-10,0.000100553)}; %1024  EOC: 1.4175
        \legend{%
$|E(\psi^N) - E(\phi^N)|$,%$ (t,\cdot)- \phi(t,\cdot) \|_{L^2(\DOmega)}$,
$\| \psi^N- u ( t^N)\|_{H^1(\DOmega)}$,
%Component \textbf{d},
}
    \end{axis}
\end{tikzpicture}
\caption{{\it Model Problem 2. Left: Consistency error $\psi^N-\phi^N$ in the $L^2$- and the $H^1$-norm for various step sizes $\tau$. Right: Comparison of the exact $H^1$-error for DS-$K5$ with the estimated error using the consistency indicator $|E(\psi^N) - E(\phi^N)|$.}}
\label{model-problem-2-consistency-errors}
\end{figure}
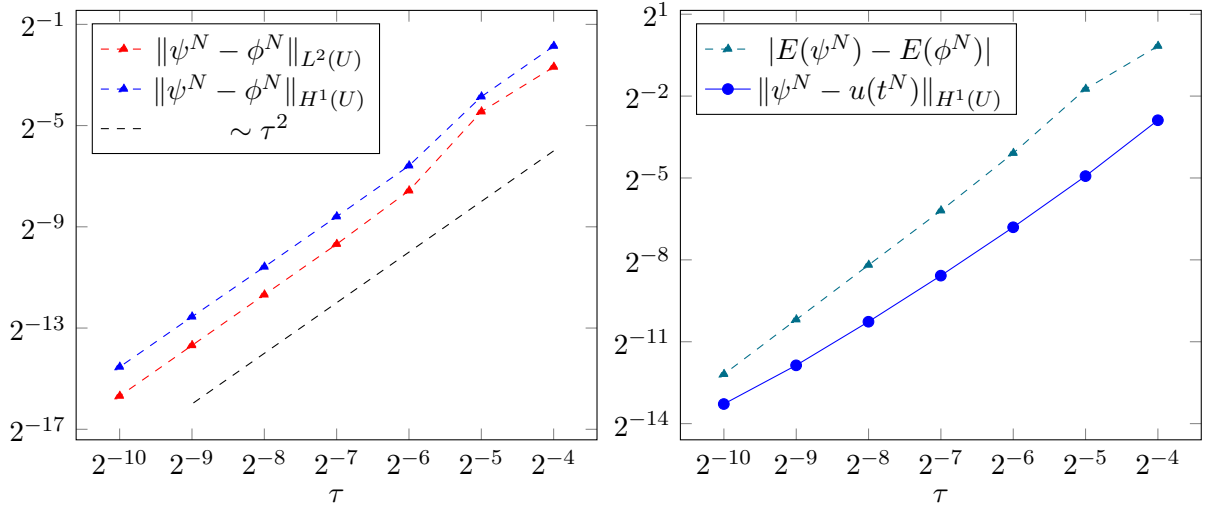

In the light of these reduced convergence orders it is interesting to investigate if the consistency error $\psi^N-\phi^N$ suffers from a similar effect. Fortunately, this does not seem to be the case as can be seen from Figure \ref{model-problem-2-consistency-errors}, where we can verify a quadratic convergence speed for both the $L^2$- and the $H^1$-consistency error. This raises the question if the previously introduced error indicator $\eta^N=|E(\psi^N) - E(\phi^N)|\sim \tau^2$ is still useful. An answer is given by the right graph in Figure \ref{model-problem-2-consistency-errors} where we observe that even though $\eta^N$ overestimates the convergence speed, it still gives a useful upper bound in the investigated regime for $\tau$. 

% Model Problem 2 - Energy and mass evolution
% Method: DS-K5

%\input{model-problem-2-energy-mass-plots}

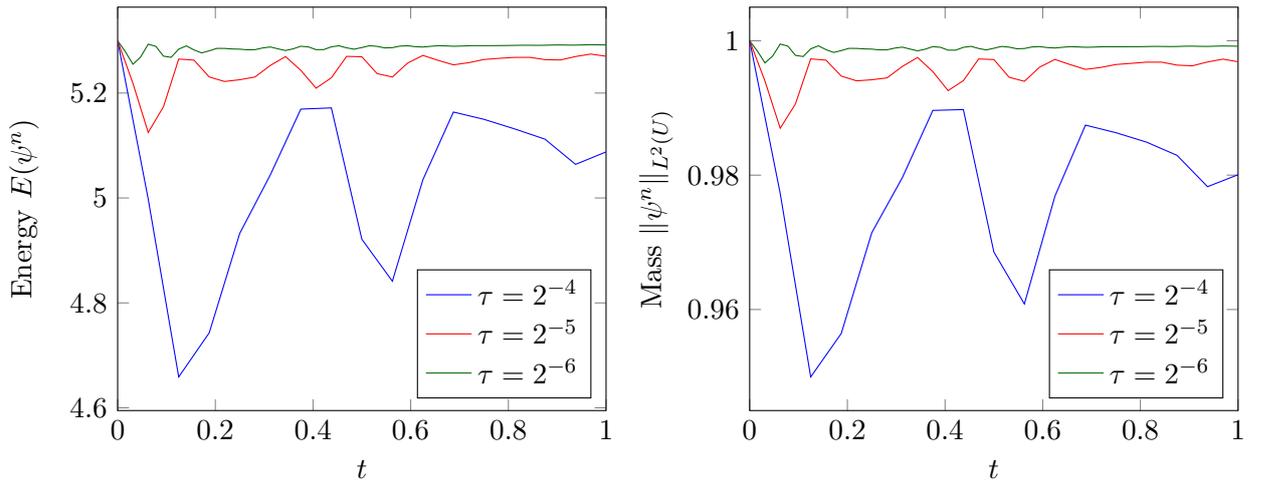
\begin{figure}
\begin{tikzpicture}
    \begin{axis}[
       legend pos=south east,
        height=0.45\textwidth,
        width=0.52\textwidth,
        xmax   = 1,  
        xmin   = 0, %ymax = 1, ymin = 0
        xlabel=$t$,
        ylabel=Energy $E(\psi^n)$]
%         \addplot[mark=*,blue] coordinates {
         \addplot[blue] coordinates {
        (0,5.29964)
        (0.0625,4.99813)
        (0.125,4.6588)
        (0.1875,4.74287)
        (0.25,4.93246)
        (0.3125,5.04382)
        (0.375,5.16935)
        (0.4375,5.1717)
        (0.5,4.92116)
        (0.5625,4.84128)
        (0.625,5.03474)
        (0.6875,5.16381)
        (0.75,5.14989)
        (0.8125,5.13173)
        (0.875,5.11203)
        (0.9375,5.06401)
        (1,5.08777)};
%         \addplot[mark=x,red] coordinates {
         \addplot[red] coordinates {
(0, 5.29964)
(0.03125, 5.21843)
(0.0625, 5.1249)
(0.09375, 5.17388)
(0.125, 5.26484)
(0.15625, 5.26273)
(0.1875, 5.23073)
(0.21875, 5.22204)
(0.25, 5.22537)
(0.28125, 5.23068)
(0.3125, 5.25254)
(0.34375, 5.26954)
(0.375, 5.2433)
(0.40625, 5.20928)
(0.4375, 5.22908)
(0.46875, 5.26994)
(0.5, 5.2691)
(0.53125, 5.23703)
(0.5625, 5.23059)
(0.59375, 5.25722)
(0.625, 5.27176)
(0.65625, 5.2623)
(0.6875, 5.25378)
(0.71875, 5.25789)
(0.75, 5.26394)
(0.78125, 5.26621)
(0.8125, 5.26799)
(0.84375, 5.26807)
(0.875, 5.26368)
(0.90625, 5.26299)
(0.9375, 5.2696)
(0.96875, 5.27433)
(1, 5.27014)};
        %
%          \addplot[mark=x,dark-green] coordinates {
          \addplot[dark-green] coordinates {
        (0, 5.29964)
(0.015625, 5.27893)
(0.03125, 5.25495)
(0.046875, 5.2686)
(0.0625, 5.29313)
(0.078125, 5.28912)
(0.09375, 5.27038)
(0.109375, 5.26805)
(0.125, 5.2837)
(0.140625, 5.28994)
(0.15625, 5.28219)
(0.171875, 5.27658)
(0.1875, 5.28043)
(0.203125, 5.28514)
(0.21875, 5.28507)
(0.234375, 5.28403)
(0.25, 5.28365)
(0.265625, 5.28242)
(0.28125, 5.28271)
(0.296875, 5.28597)
(0.3125, 5.28782)
(0.328125, 5.28432)
(0.34375, 5.28106)
(0.359375, 5.28431)
(0.375, 5.28923)
(0.390625, 5.28793)
(0.40625, 5.28273)
(0.421875, 5.28284)
(0.4375, 5.28815)
(0.453125, 5.29011)
(0.46875, 5.28618)
(0.484375, 5.28361)
(0.5, 5.28672)
(0.515625, 5.29025)
(0.53125, 5.28898)
(0.546875, 5.28603)
(0.5625, 5.28656)
(0.578125, 5.28948)
(0.59375, 5.29018)
(0.609375, 5.28844)
(0.625, 5.2877)
(0.640625, 5.28909)
(0.65625, 5.29028)
(0.671875, 5.28988)
(0.6875, 5.28921)
(0.703125, 5.28952)
(0.71875, 5.29021)
(0.734375, 5.29042)
(0.75, 5.29033)
(0.765625, 5.29038)
(0.78125, 5.29051)
(0.796875, 5.29063)
(0.8125, 5.29088)
(0.828125, 5.29114)
(0.84375, 5.29111)
(0.859375, 5.29094)
(0.875, 5.29112)
(0.890625, 5.29158)
(0.90625, 5.29171)
(0.921875, 5.29143)
(0.9375, 5.29137)
(0.953125, 5.29178)
(0.96875, 5.29212)
(0.984375, 5.29197)
(1, 5.29173)};
        \legend{%
$\tau=2^{-4}$,
$\tau=2^{-5}$,
$\tau=2^{-6}$,
%Component \textbf{d},
}
    \end{axis}
\end{tikzpicture}
\begin{tikzpicture}
    \begin{axis}[
       legend pos=south east,
        height=0.45\textwidth,
        width=0.52\textwidth,
        xmax   = 1,  
        xmin   = 0, %ymax = 1, ymin = 0
        xlabel=$t$,
        ylabel=Mass $\| \psi^n \|_{L^2(\DOmega)}$]
%         \addplot[mark=*,blue] coordinates {
         \addplot[,blue] coordinates {
        (0,1)
        (0.0625, 0.977289)
        (0.125, 0.949959)
        (0.1875, 0.956419)
        (0.25, 0.971423)
        (0.3125, 0.979618)
        (0.375, 0.989635)
        (0.4375, 0.989764)
        (0.5, 0.968556)
        (0.5625, 0.960823)
        (0.625, 0.976859)
        (0.6875, 0.987449)
        (0.75, 0.986316)
        (0.8125, 0.984915)
        (0.875, 0.982944)
        (0.9375, 0.978284)
        (1, 0.980061)};
 %         \addplot[mark=x,red] coordinates {
         \addplot[red] coordinates {
(0, 1)
(0.03125, 0.994009)
(0.0625, 0.986982)
(0.09375, 0.990599)
(0.125, 0.997305)
(0.15625, 0.997118)
(0.1875, 0.994751)
(0.21875, 0.994041)
(0.25, 0.994197)
(0.28125, 0.994482)
(0.3125, 0.996172)
(0.34375, 0.997519)
(0.375, 0.995404)
(0.40625, 0.992585)
(0.4375, 0.994066)
(0.46875, 0.997291)
(0.5, 0.997212)
(0.53125, 0.994569)
(0.5625, 0.993965)
(0.59375, 0.996074)
(0.625, 0.997226)
(0.65625, 0.996452)
(0.6875, 0.995743)
(0.71875, 0.99603)
(0.75, 0.99647)
(0.78125, 0.996639)
(0.8125, 0.996821)
(0.84375, 0.996816)
(0.875, 0.996379)
(0.90625, 0.996273)
(0.9375, 0.996847)
(0.96875, 0.997268)
(1, 0.996878)};
        %
 %          \addplot[mark=x,dark-green] coordinates {
          \addplot[dark-green] coordinates {
        (0, 1)
(0.015625, 0.998481)
(0.03125, 0.996714)
(0.046875, 0.997714)
(0.0625, 0.999512)
(0.078125, 0.999215)
(0.09375, 0.997832)
(0.109375, 0.997653)
(0.125, 0.998798)
(0.140625, 0.999251)
(0.15625, 0.99868)
(0.171875, 0.998264)
(0.1875, 0.998537)
(0.203125, 0.998871)
(0.21875, 0.998863)
(0.234375, 0.998796)
(0.25, 0.99876)
(0.265625, 0.998647)
(0.28125, 0.998655)
(0.296875, 0.99891)
(0.3125, 0.999056)
(0.328125, 0.998772)
(0.34375, 0.998502)
(0.359375, 0.998749)
(0.375, 0.999136)
(0.390625, 0.999029)
(0.40625, 0.998607)
(0.421875, 0.998605)
(0.4375, 0.999022)
(0.453125, 0.999175)
(0.46875, 0.998856)
(0.484375, 0.998642)
(0.5, 0.998884)
(0.515625, 0.999162)
(0.53125, 0.999057)
(0.546875, 0.998814)
(0.5625, 0.99885)
(0.578125, 0.99908)
(0.59375, 0.999133)
(0.609375, 0.998989)
(0.625, 0.998926)
(0.640625, 0.999033)
(0.65625, 0.999124)
(0.671875, 0.999089)
(0.6875, 0.999034)
(0.703125, 0.999055)
(0.71875, 0.999105)
(0.734375, 0.999118)
(0.75, 0.999111)
(0.765625, 0.999114)
(0.78125, 0.999118)
(0.796875, 0.999123)
(0.8125, 0.999143)
(0.828125, 0.999166)
(0.84375, 0.999159)
(0.859375, 0.999139)
(0.875, 0.999153)
(0.890625, 0.999192)
(0.90625, 0.999202)
(0.921875, 0.999173)
(0.9375, 0.999164)
(0.953125, 0.999199)
(0.96875, 0.999228)
(0.984375, 0.999212)
(1, 0.999189)};
        \legend{%
$\tau=2^{-4}$,
$\tau=2^{-5}$,
$\tau=2^{-6}$,
%Component \textbf{d},
}
    \end{axis}
\end{tikzpicture}
\caption{{\it Model Problem 2. Results obtained for DS-$K5$ and three different step sizes $\tau$. Left: Variations of the energy $E$ (cf. \eqref{energy-mod-prob-2}) over time. Right: Variations of the mass $m(\psi^n)=\| \psi^n \|_{L^2(\DOmega)}$ over time.}}
\label{model-problem-2-energy-mass-plots}
\end{figure}

Finally, we will again check the evolution of mass and energy over time. The corresponding results can be found in Figure \ref{model-problem-2-energy-mass-plots}, where we can observe that the overall picture is as in the test for Model Problem 1. Both mass and energy are oscillating, where the amplitude of the oscillations gets smaller with $\tau$. For $\tau=2^{-6}$, both energy and mass are almost constant. Again, the graphs for the mass and energy evolution look almost identical, which suggest some hidden properties of the method and maybe even the existence of an appropriate scaling that allows to simultaneously reconstruct the exact mass and energy. This remains to be investigated in more detail in future works.

\section{Summary and Conclusions}

Based on a recently developed approach for Born--Oppenheimer molecular dynamics, we have developed an 
extended ``shadow'' Lagrangian dynamics for quantum states of superfluids. The equations of motion
are given by two coupled equations, a Schr\"odinger-type equation for the quantum state and a wave equation 
for an extended auxiliary field variable. The equations are coupled in a nonlinear way, but each equation 
individually is linear with respect to the variable that it defines. Thanks to this linearization
the system can be easily discretized using linear time stepping methods, where we used a Crank-Nicolson-type 
approach for the Schr\"odinger equation and an extended leapfrog scheme for the wave equation. This allows
for a fast and stable integration of the equations of motion that is competitive to current state-of-the-art methods.
We also found that the difference between the quantum state and the extended field variable defines a consistency error
that can be used to estimate the numerical error on-the-fly without any additional costs.

%\section{Acknowledgements}
%
%
%AMNN is most grateful for the hospitality and pleasant environment at the Division of Scientific Computing at
%the department of Information Technology at Uppsala University, where he stayed during his initial participation of this work.
%This work is supported by the U.S. Department of Energy Office of Basic Energy Sciences (FWP LANLE8AN)
%and by the U.S. Department of Energy through the Los Alamos National Laboratory.
%Los Alamos National Laboratory is operated by Triad National Security, LLC, for the National Nuclear Security
%Administration of the U.S. Department of Energy Contract No. 892333218NCA000001.

\def\cprime{$'$}

\end{document}